\documentclass[a4paper, 11pt, reqno]{amsart}

\usepackage[margin={25mm, 25mm}]{geometry}
\usepackage{amsmath}
\usepackage{amssymb}
\usepackage{amsthm}
\usepackage{float}
\usepackage{mathrsfs}
\usepackage{url}
\usepackage{graphicx}
\usepackage{color}
\usepackage[legacycolonsymbols]{mathtools} 
\usepackage{tikz}
\usepackage{bbm}
\usepackage{comment}
\usepackage{bm}
\usepackage{multirow}
\usepackage{mathdots}
\usepackage[colorlinks=true, linkcolor=blue, citecolor=blue, urlcolor=blue, breaklinks=true]{hyperref}

\theoremstyle{definition}
\newtheorem{defi}{Definition}[section]
\newtheorem{ex}[defi]{Example}

\theoremstyle{plain}
\newtheorem{thm}[defi]{Theorem}
\newtheorem{prop}[defi]{Proposition}
\newtheorem{lem}[defi]{Lemma}
\newtheorem{cor}[defi]{Corollary}

\newtheorem*{claim}{Claim}

\numberwithin{equation}{section}

\newcommand{\Z}{\mathbb{Z}}
\newcommand{\Q}{\mathbb{Q}}

\newcommand{\T}{\mathbb{T}}
\newcommand{\PP}{\mathbb{P}}

\newcommand{\CF}{\mathcal{F}}

\newcommand{\CP}{\mathcal{P}}

\newcommand{\vv}{\vspace{1em}}
\newcommand{\hh}{\hspace{0.5em}}

\newcommand{\bx}{{\mathbf{x}}}
\newcommand{\by}{{\mathbf{y}}}
\newcommand{\bz}{\mathbf{z}}
\newcommand{\bp}{{\mathbf{p}}}

\newcommand{\be}{{\mathbf{e}}}
\newcommand{\bb}{{\mathbf{b}}}

\newcommand{\bg}{{\mathbf{g}}}
\newcommand{\bc}{{\mathbf{c}}}
\newcommand{\br}{\mathbf{r}}
\newcommand{\bu}{\mathbf{u}}

\newcommand{\bze}{\boldsymbol{0}}
\newcommand{\bsig}{\mathbf{\Sigma}}

\newcommand{\ee}{\end{equation*}}

\newcommand{\diag}{{\mathrm{diag}}}

\newcommand{\op}{\oplus}

\newcommand{\Tp}{\mathrm{Trop}}
\newcommand{\one}{\mathbbm{1}}
\newcommand{\hy}{\hat{y}}
\newcommand{\hby}{\hat{\by}}
\newcommand{\vep}{\varepsilon}

\allowdisplaybreaks

\title[SEPARATION FORMULAS FOR NON-NORMALIZED GENERALIZED CLUSTER ALGEBRAS]{Separation formulas for non-normalized generalized cluster algebras}
\author{Koji Amagai}
\address{Graduate School of Mathematics, Nagoya University, Chikusa-ku, Nagoya, 464-8604, Japan}
\email{koji.amagai.h5@math.nagoya-u.ac.jp}
\date{}

\begin{document}

\begin{abstract}
We study a class of not necessarily normalized generalized cluster patterns which are uniquely determined by the initial seed and mutations. For such generalized cluster patterns, we provide the separation formulas for cluster variables and coefficients.
\end{abstract}

\maketitle

\tableofcontents

\section{Introduction}
Cluster algebras are a class of algebras introduced by Fomin and Zelevinsky \cite{CA1}. Subsequently, Chekhov and Shapiro \cite{CS14} introduced generalized cluster algebras, which are a natural generalization of ordinary cluster algebras, by extending the exchange relations from a binomial to a polynomial. Moreover, it is known that fundamental properties established in ordinary cluster algebras such as the Laurent phenomenon, the finite type classification \cite{CS14}, the separation formulas \cite{Nak15,NR16}, and the Laurent positivity \cite{Mou24,BLM25a,BLM25b} hold in generalized cluster algebras. The greedy basis of a rank 2 cluster algebra is extended to a rank 2 generalized cluster algebra \cite{Rup15}. The structure of generalized cluster algebras appears in various contexts, such as exact WKB analysis \cite{IN16}, the Drinfeld double of $GL_n$ \cite{GSV18, GSV22a, GSV22b}, and the Badlands in (ordinary) cluster scattering diagrams of acyclic affine type \cite{RS26a,RS26b}. More recently, Ramos and Whiting \cite{RW25} showed that a generalized cluster algebra of geometric type is isomorphic to a subquotient of an ordinary cluster algebra, and Akagi and Nakanishi \cite{AN26} extended this result to generalized cluster algebras with arbitrary coefficients. 

The separation formulas introduced in \cite{CA4} are important for ordinary cluster algebras and their underlying ordinary cluster patterns. In particular, the $C$-, $G$-matrices, and the $F$-polynomials appearing in these formulas are major objects for understanding the structure of cluster algebras. For instance, the periodicity of cluster variables coincides with that of $G$-matrices \cite{Nak21}. Nakanishi and Rupel established the separation formulas for generalized cluster patterns \cite{Nak15,NR16}. However, the above separation formulas only hold for \textit{normalized} (ordinary or generalized) cluster patterns, where the normalization condition is satisfied by their coefficients. On the other hand, non-normalized (ordinary or generalized) cluster algebras are also important objects to study (e.g., \cite{FT18,MOT26} for ordinary cluster algebras, \cite{GSV18,GSV22a,GSV22b} for generalized ones).

In this paper, as the main results, we extend the separation formulas for normalized (ordinary or generalized) cluster patterns to non-normalized ones (Theorems \ref{prop:non-sep} and \ref{GCA-sep-p}). For this purpose, we need to introduce several new notions and formulations. Firstly, in contrast to the normalized (ordinary or generalized) cluster patterns, non-normalized (ordinary or generalized) cluster patterns are not uniquely determined by their initial seed and the mutations. Therefore, we first extend the usual normalization condition to the \textit{$\CP$-normalization condition} (Definitions \ref{P-norm cond} and \ref{GCA_P-norm cond}). Then, an (ordinary or generalized) cluster pattern satisfying this condition, which we call a \textit{$\CP$-normalized (ordinary or generalized) cluster pattern}, is uniquely determined by their initial seed and the mutations and the separation formulations are available. Secondly, for generalized cluster patterns, we introduce $y$-variables whose formulation is different from that in \cite{Nak15, NR16}. This construction allows us to define $y$-variables without assuming the existence of their roots in \cite{Nak15,NR16}, as detailed in Section \ref{corres}. Finally, in the resulting separation formula of a cluster variable, the numerator is given by the $F$-polynomials as usual, while the denominator is given by a new function (a \textit{$J$-function}), which is no longer a polynomial.

The organization of this paper is as follows. First, in Section \ref{reform}, we concentrate on \textit{ordinary} cluster patterns. We introduce the $\CP$-normalization condition. Then, we present the separation formulas for $\CP$-normalized cluster patterns (Theorem \ref{prop:non-sep}), which generalize the formulas in \cite{CA4}.

In Section \ref{P-GCP}, we give the basic notions of the generalized cluster patterns. We introduce the $\CP$-normalization condition as in the ordinary case, and introduce $y$-variables for $\CP$-normalized generalized cluster patterns.

In Section \ref{GCA-Sep}, in parallel to the ordinary case \cite{CA4}, we introduce generalized cluster patterns with principal coefficients, and construct the $C$-, $G$-matrices, the $F$-polynomials, and the $J$-functions from the information of the $x$- and $y$-variables defined in Section \ref{P-GCP}. Furthermore, after analyzing the detailed block structure of the $C$-matrices, we give the separation formulas for generalized cluster patterns (Theorems \ref{prop:GCA-sep} and \ref{GCA-sep-p}).

Finally, in Section \ref{Further}, we show several further properties. In Section \ref{subsec:corres}, we investigate the correspondence of the $C$-, $G$-matrices, and the $F$-polynomials between a generalized cluster pattern and an associated ordinary one. We also show that dualities of the same form as those known for ordinary cluster patterns appear in generalized cluster patterns as well in Section \ref{subsec:dual}. In Section \ref{corres}, we explain the correspondence between our formulation and the one in \cite{Nak15,NR16}.

\section*{Acknowledgments}
The author is grateful to Professor Tomoki Nakanishi for discussions, careful reading, and advice.

\section{Separation formulas for \texorpdfstring{$\CP$}{P}-normalized cluster algebras} \label{reform}

In this section, we recall basic notions for cluster patterns following \cite{CA2} and \cite{CA4}. Under the notions, we introduce the $P$-patterns and $\CP$-normalized cluster patterns. Furthermore, we define the $J$-functions, and using them, we provide the separation formulas for $\CP$-normalized cluster patterns.

Although ordinary cluster algebras are a special case of generalized cluster algebras, they remain important; therefore, we discuss ordinary cluster algebras as follows.

\subsection{Non-normalized cluster patterns}\label{subsec:non-clus-pat}
We recall the definition of non-normalized cluster patterns and their mutations based on \cite{CA2} and \cite{CA4}. \vv

\begin{defi}[Semifields]
Let $\PP$ be a multiplicative abelian group with respect to the product $\cdot$, endowed with an addition $\op$ satisfying the following conditions. Then, the triplet $(\PP,\op,\cdot)$ is called a \emph{semifield}. 
\begin{enumerate}
    \item $a\op b=b\op a$.
    \item $a(b\op c)=ab\op ac$.
    \item $a\op(b\op c)=(a\op b)\op c$.
\end{enumerate}
\end{defi}

The group ring $\Z\PP$ of $\PP$ is an integral domain \cite{CA1}. Thus, its field of the fractions $\Q\PP$ is well-defined. Let $\CF$ be a field that is isomorphic to the field of the rational functions in $n$ variables with coefficients in $\Q\PP$, which we call the \emph{ambient field}.

An $n \times n$ integer matrix $B=(b_{ij})$ is called \emph{skew-symmetrizable} if there exists a diagonal matrix $D$ with positive rational diagonal entries such that $DB=-(DB)^T$. The matrix $D$ is called a \emph{skew-symmetrizer}. In particular, $b_{ii}=0$ if $B$ is a skew-symmetrizable matrix.

\begin{defi}[Seeds] A \emph{seed} is a triplet $\Sigma^p=(\bx,\bp,B)$ given by the following:
\begin{enumerate}
    \item $\bx=(x_1,\dots,x_n)$ is an $n$-tuple of algebraically independent elements over $\Q\PP$ generating $\CF$, whose elements are called \emph{cluster variables} or \emph{$x$-variables}.
    \item $\bp=(p_1^+,p_1^-,\dots,p_n^+,p_n^-)$ is a $2n$-tuple of elements in $\PP$ called a \emph{coefficient tuple}, and its components are called \emph{$p$-variables}.
    \item $B$ is an $n \times n$ skew-symmetrizable matrix, called an \emph{exchange matrix}.
\end{enumerate}
\end{defi}

Here, we use the notation $\Sigma^p$ to emphasize that coefficients are given by $p$-variables.

\begin{defi}[Seed mutations] A seed $(\Sigma^p)'=(\bx',\bp',B')$ is called a mutation of $\Sigma^p$ in direction $k=1,\dots,n$ if it satisfies the following rules:
\begin{align}
&b_{ij}'=\begin{cases}
-b_{ij} & \text{$(i=k$ or $j=k)$},\\
b_{ij}+b_{ik}[b_{kj}]_++[-b_{ik}]_+b_{kj} & \text{(otherwise)},
\end{cases}\label{B-mut}\\
&(p_{k}^\pm)'=p_{k}^\mp\label{p-mut1},\\
&\dfrac{(p_{i}^+)'}{(p_{i}^-)'}=\dfrac{p_{i}^+}{p_{i}^-}(p_{k}^-)^{-[-b_{ki}]_+}(p_{k}^+)^{[b_{ki}]_+} \quad \text{$(i\neq k)$},\label{p-mut2}\\
&x_{i}'=\begin{cases}
\displaystyle x_{k}^{-1}\left(p_{k}^+\left(\prod_{j=1}^nx_{j}^{[b_{jk}]_+}\right)+p_{k}^-\left(\prod_{j=1}^nx_{j}^{[-b_{jk}]_+}\right)\right) & \text{$(i=k)$},\\
x_{i} & \text{$(i\neq k)$},
\end{cases}\label{exchange-rel}
\end{align}
where $[a]_+=\max\{a,0\}$.
\end{defi}

It is important that $\bp'$ is not uniquely determined from $\bp$ by \eqref{p-mut1} and \eqref{p-mut2} because they only specify the ratio $(p_i^+)'/(p_i^-)'$ for $i\neq k$.

Noting that $a=[a]_+-[-a]_+$, we obtain
\begin{align}
b_{ik}[b_{kj}]_+&=b_{ik}(b_{kj}+[-b_{kj}]_+),\label{Bcal1}\\
[-b_{ik}]_+b_{kj}&=([b_{ik}]_+-b_{ik})b_{kj}.\label{Bcal2}
\end{align}
From \eqref{Bcal1} and \eqref{Bcal2}, it follows that
\begin{align*}
b_{ik}[b_{kj}]_++[-b_{ik}]_+b_{kj}&=b_{ik}(b_{kj}+[-b_{kj}]_+)+([b_{ik}]_+-b_{ik})b_{kj}\\
&=b_{ik}[-b_{kj}]_++[b_{ik}]_+b_{kj}\\
&=-b_{ik}'[b_{kj}']_+-[-b_{ik}']_+b_{kj}'.
\end{align*}
Combining it with \eqref{B-mut} yields
\begin{align}
b_{ij}=\begin{cases}
-b_{ij}' & \text{$(i=k$ or $j=k)$},\\
b_{ij}'+b_{ik}'[b_{kj}']_++[-b_{ik}']_+b_{kj}' & \text{(otherwise)}.
\end{cases}
\end{align}
Furthermore, \eqref{p-mut1} and the first case in \eqref{B-mut} imply that \eqref{p-mut2} is equivalent to:
\begin{align}
\dfrac{p_{i}^+}{p_{i}^-}=\dfrac{(p_{i}^+)'}{(p_{i}^-)'}((p_{k}^-)')^{-[-b_{ki}']_+}((p_{k}^+)')^{[b_{ki}']_+} \quad \text{$(i\neq k)$}.
\end{align}
Moreover, noting that $b_{kk}=0$, by \eqref{B-mut}, \eqref{p-mut1}, and the second case of \eqref{exchange-rel}, the first case of \eqref{exchange-rel} is equivalent to:
\begin{align}
x_k=\displaystyle (x_{k}')^{-1}\left((p_{k}^-)'\left(\prod_{j=1}^n(x_{j}')^{[-b_{jk}']_+}\right)+(p_{k}^+)'\left(\prod_{j=1}^n(x_{j}')^{[b_{jk}']_+}\right)\right).
\end{align}
Therefore, if $(\Sigma^p)'$ is a mutation of $\Sigma^p$ in direction $k$, then $\Sigma^p$ is also a mutation of $(\Sigma^p)'$ in direction $k$.

Let $\T_n$ be an \emph{$n$-regular tree} whose edges emanating from each vertex are uniquely labeled by $1, \dots, n$. By abuse of notation, we also use $\T_n$ to denote its set of the vertices. Two vertices $t, t' \in \T_n$ are called \emph{$k$-adjacent} if they are connected by an edge labeled by $k$.

\begin{defi}[Non-normalized cluster patterns] 
Suppose that we attach a seed $\Sigma^p_t=(\bx_t,\bp_t,B_t)$ to each vertex $t \in \mathbb{T}_n$. The family of the seeds $\bsig^p=(\Sigma^p_t)_{t\in\T_n}$ is called a \emph{(non-normalized) cluster pattern} if for any two vertices $t$ and $t'$ connected by an edge labeled $k$, the seed $\Sigma^p_{t'}$ is the mutation of $\Sigma^p_{t}$ in direction $k$.
\end{defi}

As previously noted, the mutations do not uniquely determine the subsequent $p$-variables. Consequently, even if we fix an initial seed $\Sigma_{t_0}^p$ at an initial vertex $t_0\in\T_n$, a non-normalized cluster pattern is not uniquely determined by the seed.

\subsection{\texorpdfstring{$\CP$}{P}-normalized cluster patterns}\label{subsec:P-pet}
As mentioned in Section \ref{subsec:non-clus-pat}, in general, a non-normalized cluster pattern is not uniquely determined by an initial seed. In this subsection, we introduce a $P$-pattern and define a $\CP$-normalized cluster pattern as a non-normalized cluster pattern uniquely determined by an initial seed.

\begin{defi}\label{P-pattern} Let $\CP=\{p_{i;t}\in\PP \mid i\in\{1,\dots,n\},t\in\T_n\}$ be a family of elements such that $p_{k;t}=p_{k;t'}$ for any $k$-adjacent vertices $t,t'\in\T_n$. We call $\CP$ a \emph{$P$-pattern}. When $p_{i;t}=1$ for all $i\in\{1,\dots,n\}$ and $t\in\T_n$, we call $\CP$ the \emph{trivial $P$-pattern} and denote it by $\CP_\mathrm{triv}$. 
\end{defi}

\begin{defi}\label{P-norm cond} Let $\CP$ be a $P$-pattern. Suppose that a cluster pattern $\bsig^p=(\Sigma_t^p)_{t\in\T_n}$ satisfies the following condition, called the \emph{$\CP$-normalization condition}:
\begin{align}
p_{i;t}^+\op p_{i;t}^-=p_{i;t},\label{P-nor}
\end{align}
for any $i\in\{1,\dots,n\}$ and $t\in\T_n$. In this case, we call the cluster pattern $\bsig^p$ a \emph{$\CP$-normalized cluster pattern}. When $\CP$ is $\CP_{\mathrm{triv}}$, the condition \eqref{P-nor} is the (standard) normalization condition in \cite{CA2}. Also, a $\CP_{\mathrm{triv}}$-normalized cluster pattern is simply called a normalized cluster pattern.
\end{defi}

For a $\CP$-normalized cluster pattern $\bsig^p$, let $y_{i;t}=p_{i;t}^+/p_{i;t}^-$ which is called a \emph{$y$-variable}. By the definition, we have $p_{i;t}^+=p_{i;t}^-y_{i;t}$ for any $i=1,\dots,n$ and $t\in\T_n$. Combining it with \eqref{P-nor}, we obtain $p_{i;t}^-(1\op y_{i;t})=p_{i;t}$, and then, we have the following relation:
\begin{align}
p_{i;t}^-=\dfrac{p_{i;t}}{1\op y_{i;t}}.\label{p to y1}
\end{align}
Moreover, by \eqref{P-nor} and \eqref{p to y1}, we obtain
\begin{align}
p_{i;t}^+=\dfrac{p_{i;t}y_{i;t}}{1\op y_{i;t}}. \label{p to y2}
\end{align}
\begin{prop}\label{prop:y-mut} Suppose that $\bsig^p=(\Sigma^p_t)_{t\in\T_n}$ is a $\CP$-normalized cluster pattern and let $t'$ be $k$-adjacent to $t$. Then, we obtain the following.

\textup{(1)} Each $y$-variable $y_{i;t}$ in $\bsig^p$ mutates in direction $k$ as follows:
\begin{align}
y_{i;t'}=\begin{cases}
y_{k;t}^{-1} & \text{$(i=k)$},\\
p_{k;t}^{b_{ki;t}}y_{i;t}y_{k;t}^{[b_{ki;t}]_+}(1\op y_{k;t})^{-b_{ki;t}} & \text{$(i\neq k)$}.
\end{cases}\label{y-mut}
\end{align}

\textup{(2)} The mutation of a seed $\Sigma^p_t$ in direction $k$ is uniquely determined.
\end{prop}
\begin{proof}
(1) The first case of \eqref{y-mut} follows immediately from \eqref{p-mut1}. Next, we consider the second case of \eqref{y-mut}. 
Using \eqref{p-mut2}, \eqref{p to y1}, and \eqref{p to y2}, we calculate $y_{i;t'}$ for $i\neq k$ as follows:
\begin{align*}
y_{i;t'}&=y_{i;t}\left(\dfrac{p_{k;t}}{1\op y_{k;t}}\right)^{-[-b_{ki;t}]_+}\left(\dfrac{p_{k;t}y_{k;t}}{1\op y_{k;t}}\right)^{[b_{ki;t}]_+}\\
&=y_{i;t}y_{k;t}^{[b_{ki;t}]_+}\left(\dfrac{p_{k;t}}{1\op y_{k;t}}\right)^{[b_{ki;t}]_+-[-b_{ki;t}]_+}\\
&=p_{k;t}^{b_{ki;t}}y_{i;t}y_{k;t}^{[b_{ki;t}]_+}(1\op y_{k;t})^{-b_{ki;t}}.
\end{align*}
(2) By the definition of the seed mutations, it suffices to show that the mutation of the $p$-variables in direction $k$ is uniquely determined. First, it is clear that $p_{k;t'}^\pm$ are uniquely determined by \eqref{p-mut1}. Next, we show that $p_{i;t'}^\pm$ are uniquely determined for $i\neq k$. By (1), the mutation of the $n$-tuple of $y$-variables $\by_t$ in direction $k$ is uniquely determined. Now, the $p$-variables $p_{1;t'}^\pm,\dots,p_{n;t'}^\pm$ are expressed in terms of each $p_{i;t'}\in\CP$ and the $y$-variables $y_{1;t'},\dots,y_{n;t'}$ as in \eqref{p to y1} and \eqref{p to y2}. Since each $p_{i;t'}$ is fixed, the claim (2) holds.
\end{proof}

\begin{prop} A $\CP$-normalized cluster pattern is uniquely determined by the $P$-pattern $\CP$ and an initial seed $\Sigma_{t_0}^p$. 
\end{prop}
\begin{proof}
It immediately follows from Proposition \ref{prop:y-mut} (2).
\end{proof}

Since the $p$-variables can be recovered from the $y$-variables by \eqref{p to y1} and \eqref{p to y2}, it suffices to consider the $y$-variables instead of the $p$-variables. Accordingly, we examine the seed $\Sigma_t=\Sigma_t^y=(\bx_t,\by_t,B_t)$ and the $\CP$-normalized cluster pattern $\bsig=(\Sigma_t)_{t\in\T_n}$, rather than their $p$-counterparts $\Sigma_t^p$ and $\bsig^p$.

Now, we set $\hy_{i;t}=y_{i;t}\prod_{j=1}^nx_{j;t}^{b_{ji;t}}$, which is called a \emph{$\hy$-variable}, and denote the $n$-tuple of these variables by $\mathbf{\hy}_t=(\hy_{1;t},\dots,\hy_{n;t})$. Then, the first case of \eqref{exchange-rel} is rewritten as follows:
\begin{align}
x_{k;t'}=p_{k;t}x_{k;t}^{-1}\left(\prod_{j=1}^nx_{j;t}^{[-b_{jk;t}]_+}\right)\frac{1+\hy_{k;t}}{1\op y_{k;t}}\label{exchange-rel2}.
\end{align}
Furthermore, $\hy$-variables mutate in direction $k$ as follows:

\begin{prop}[{\cite[Proposition 2.9]{FT18}}] Let $t'$ be $k$-adjacent to $t$. Then, the following mutation formula holds.
\begin{align}
\hy_{i;t'}=\begin{cases}
\hy_{k;t}^{-1} & \text{$(i=k)$},\\
\hy_{i;t}\hy_{k;t}^{[b_{ki;t}]_+}(1+ \hy_{k;t})^{-b_{ki;t}} & \text{$(i\neq k)$}.
\end{cases}\label{hy-mut}
\end{align}
\end{prop}
\begin{proof}
By the second case of \eqref{exchange-rel}, we have
\begin{align*}
\hy_{i;t'}=y_{i;t'}\left(\prod_{\substack{j=1\\j\neq k}}^nx_{j;t}^{b_{ji;t'}}\right)x_{k;t'}^{b_{ki;t'}}.
\end{align*}
First, $\hy_{k;t'}=\hy_{k;t}^{-1}$ holds since $y_{k;t'}=y_{k;t}^{-1}$, $b_{jk;t'}=-b_{jk;t}$, and $b_{kk;t}=0$. Next, we show the second case of \eqref{hy-mut}. For $i\neq k$, by the mutations of the $x$-, $y$-variables, and the exchange matrix \eqref{exchange-rel2}, \eqref{y-mut}, \eqref{B-mut}, we have:
\begin{align*}
\hy_{i;t'}&=p_{k;t}^{b_{ki;t}}y_{i;t}y_{k;t}^{[b_{ki;t}]_+}(1\op y_{k;t})^{-b_{ki;t}}\\
&\quad \times \prod_{\substack{j=1\\j\neq k}}^nx_{j;t}^{b_{ji;t}+b_{jk;t}[b_{ki;t}]_++[-b_{jk;t}]_+b_{ki;t}}\left\{p_{k;t}x_{k;t}^{-1}\left(\prod_{j=1}^nx_{j;t}^{[-b_{jk;t}]_+}\right)\frac{1+\hy_{k;t}}{1\op y_{k;t}}\right\}^{-b_{ki;t}}\\
&=y_{i;t}\left(\prod_{j=1}^nx_{j;t}^{b_{ji;t}}\right) \times \left(y_{k;t}^{[b_{ki;t}]_+}\prod_{\substack{j=1\\j\neq k}}^nx_{j;t}^{b_{jk;t}[b_{ki;t}]_+}\right)\times (1+ \hy_{k;t})^{-b_{ki;t}}\\
&=\hy_{i;t}\hy_{k;t}^{[b_{ki;t}]_+}(1+ \hy_{k;t})^{-b_{ki;t}}.
\end{align*}
This completes the proof.
\end{proof}

We note that the mutation of $\hy$-variables does not depend on $\CP$. Also, when $\CP$ is the trivial $P$-pattern, the form of the mutation formula for $\hy$-variables coincides with that for $y$-variables.

\subsection{Separation formulas for \texorpdfstring{$\CP$}{P}-normalized cluster patterns}

In the following, we extend the separation formulas for normalized cluster patterns given in \cite{CA4} and provide the separation formulas for $\CP$-normalized cluster patterns. First, we define the $C$-, $G$-matrices, the $F$-polynomials, and the $J$-functions associated with a $\CP$-normalized cluster pattern $\bsig=(\Sigma_t)_{t\in\T_n}$. The $C$-, $G$-matrices, and $F$-polynomials are the same as the standard ones. Here, we write $\bx=\bx_{t_0}=(x_1,\dots,x_n)$, $\by=\by_{t_0}=(y_1,\dots,y_n)$, and $\hby=\hby_{t_0}=(\hy_1,\dots,\hy_n)$ for the initial variables at $t_0\in\T_n$.

\begin{defi}
Let $\bsig=(\Sigma_t)_{t\in\T_n}$ be a $\CP$-normalized cluster pattern with an initial vertex $t_0$, and suppose that $t'$ is $k$-adjacent to $t$. We define the \emph{$C$-, $G$-matrices}, the \emph{$F$-polynomials}, and the \emph{$J$-functions} by the following recurrence relations:
\begin{align}
&C_{t_0}=I_n ,\quad c_{ij;t'}=\begin{cases}
-c_{ik;t} & \text{$(j=k)$},\\
c_{ij;t}+c_{ik;t}[b_{kj;t}]_++[-c_{ik;t}]_+b_{kj;t} & \text{$(j\neq k)$},
\end{cases}\label{OCA-C-mut}\\
&G_{t_0}=I_n ,\quad g_{ij;t'}=\begin{cases}
-g_{ik;t}+\sum_{\ell=1}^ng_{i\ell;t}[-b_{\ell k;t}]_+-\sum_{\ell=1}^nb_{i\ell;t_0}[-c_{\ell k;t}]_+ & \text{$(j=k)$},\\
g_{ij;t} & \text{$(j\neq k)$},
\end{cases}\label{OCA-G-mut}\\
&F_{i;t_0}(\by)=1,\quad F_{i;t'}(\by)=\begin{cases}
\dfrac{M_{k;t}(\by)}{F_{k;t}(\by)} & \text{$(i=k)$},\\
F_{i;t}(\by) & \text{$(i\neq k)$},
\end{cases}\label{OCA-F-mut}\\
&J_{i;t_0}(\by,\CP)=1,\quad J_{i;t'}(\by,\CP)=\begin{cases}
\dfrac{N_{k;t}(\by,\CP)}{p_{k;t}J_{k;t}(\by,\CP)} & \text{$(i=k)$},\\
J_{i;t}(\by,\CP) & \text{$(i\neq k)$},
\end{cases}\label{OCA-J-mut}
\end{align}
where $M_{k;t}(\by)$ and $N_{k;t}(\by,\CP)$ are given by
\begin{align*}
M_{k;t}(\by)
&=\left(\prod_{j=1}^ny_j^{[-c_{jk;t}]_+}\right)\left(\prod_{j=1}^nF_{j;t}(\by)^{[-b_{jk;t}]_+}\right)
\left\{1+\left(\prod_{j=1}^ny_j^{c_{jk;t}}\right)\left(\prod_{j=1}^nF_{j;t}(\by)^{b_{jk;t}}\right)\right\},\\
N_{k;t}(\by,\CP)
&=\left(\prod_{j=1}^ny_j^{[-c_{jk;t}]_+}\right)\left(\prod_{j=1}^nJ_{j;t}(\by,\CP)^{[-b_{jk;t}]_+}\right)
\left\{1+\left(\prod_{j=1}^ny_j^{c_{jk;t}}\right)\left(\prod_{j=1}^nJ_{j;t}(\by,\CP)^{b_{jk;t}}\right)\right\}.
\end{align*}

The transformations $C_t \mapsto C_{t'}$, $G_t \mapsto G_{t'}$, $F_{i;t}(\by) \mapsto F_{i;t'}(\by)$, and $J_{i;t}(\by,\CP) \mapsto J_{i;t'}(\by,\CP)$ are called the \emph{mutations} of the $C$-, $G$-matrices, the $F$-polynomials, and the $J$-functions, respectively. The resulting families $\{C_t\}_{t\in\T_n}$, $\{G_t\}_{t\in\T_n}$, $\{F_{i;t}\}_{1\le i \le n, t\in\T_n}$, and $\{J_{i;t}\}_{1\le i \le n, t\in\T_n}$ are called the \emph{$C$-}, \emph{$G$-}, \emph{$F$-}, and \emph{$J$-patterns}, respectively.
\end{defi}

$M_{k;t}(\by)$ is obtained from $N_{k;t}(\by,\CP)$ by the specialization setting all $p_{j;t}$ to $1$: $M_{k;t}(\by)=N_{k;t}(\by,\one)$. Thus, the $F$-polynomials are obtained from the $J$-functions by specializing all $p_{j;t}$ to $1$: $F_{i;t}(\by)=J_{i;t}(\by,\one)$. Unlike the $F$-polynomials, the $J$-functions are rational functions in $\by$ and $\CP$.

In the following, we generalize the separation formulas in \cite{CA4}.

\begin{thm}[Separation formulas for $\CP$-normalized cluster patterns]\label{prop:non-sep} The following formulas hold (formula \eqref{hy-sep} is due to \textup{\cite[Proposition 3.13]{CA4}}):
\begin{align}
y_{i;t}&=\left(\prod_{j=1}^ny_j^{c_{ji;t}}\right)\prod_{j=1}^n J_{j;t}|_\PP(\by,\CP)^{b_{ji;t}},\label{y-sep}\\
\hy_{i;t}&=\left(\prod_{j=1}^n\hy_j^{c_{ji;t}}\right)\prod_{j=1}^n F_{j;t}(\hby)^{b_{ji;t}},\label{hy-sep}\\
x_{i;t}&=\left(\prod_{j=1}^nx_j^{g_{ji;t}}\right)\dfrac{F_{i;t}(\hby)}{J_{i;t}|_\PP(\by,\CP)}.\label{x-sep}
\end{align}
\end{thm}
\begin{proof}
We can prove them by a similar argument to \cite[proof of Theorem 4.16]{Nak23}, as follows.
For $t=t_0$, they clearly hold from the initial conditions of $C$-, $G$-matrices, $F$-polynomials and $J$-functions. Suppose that the claim holds at $t\in\T_n$ and let $t'$ be $k$-adjacent to $t$. It is sufficient to recover the mutation formulas of the $y$-, $\hy$-, and the $x$-variables by the mutation formulas of the $C$-,$G$-matrices, the $F$-polynomials, and the $J$-functions. First, we show the formula \eqref{y-sep}. For $i=k$, we have $c_{jk;t'}=-c_{jk;t}$, $J_{j;t'}|_\PP(\by,\CP)=J_{j;t}|_\PP(\by,\CP)$ for $j\neq k$, $b_{jk;t'}=-b_{jk;t}$ and $b_{kk;t'}=0$. Thus the formula \eqref{y-sep} holds for $i=k$. For $i\neq k$, writing the right-hand side of \eqref{y-sep} as $y_{i;t'}$, we have the following:
\begin{align*}
y_{i;t'}&=\left(\prod_{j=1}^ny_j^{c_{ji;t'}}\right)\prod_{j=1}^n J_{j;t'}|_\PP(\by,\CP)^{b_{ji;t'}}\\
&=\left(\prod_{j=1}^ny_j^{c_{ji;t}+c_{jk;t}[b_{ki;t}]_++[-c_{jk;t}]_+b_{ki;t}}\right) \times \prod_{\substack{j=1\\j\neq k}}^nJ_{j;t}|_\PP(\by,\CP)^{b_{ji;t}+b_{jk;t}[b_{ki;t}]_++[-b_{jk;t}]_+b_{ki;t}}\\
&\quad\times \left\{p_{k;t}^{-1}J_{k;t}|_\PP(\by,\CP)^{-1}\left(\displaystyle\prod_{j=1}^ny_j^{[-c_{jk;t}]_+}\right)\displaystyle\prod_{j=1}^nJ_{j;t}|_\PP(\by,\CP)^{[-b_{jk;t}]_+}\right.\\
&\quad\left.\times\left(1\op\displaystyle\left(\prod_{j=1}^ny_j^{c_{jk;t}}\right)\displaystyle\prod_{j=1}^nJ_{j;t}|_\PP(\by,\CP)^{b_{jk;t}}\right)\right\}^{-b_{ki;t}}\\
&=p_{k;t}^{b_{ki;t}}\left(\prod_{j=1}^ny_j^{c_{ji;t}}\right)J_{k;t}|_\PP(\by,\CP)^{b_{ki;t}}\prod_{\substack{j=1\\j\neq k}}^nJ_{j;t}|_\PP(\by,\CP)^{b_{ji;t}}\\
&\quad\times \left(\prod_{j=1}^ny_j^{c_{jk;t}[b_{ki;t}]_+}\right)\prod_{\substack{j=1\\j\neq k}}^nJ_{j;t}|_\PP(\by,\CP)^{b_{jk;t}[b_{ki;t}]_+}
\times\left(1\op\displaystyle\left(\prod_{j=1}^ny_j^{c_{jk;t}}\right)\displaystyle\prod_{j=1}^nJ_{j;t}|_\PP(\by,\CP)^{b_{jk;t}}\right)^{-b_{ki;t}}\\
&\overset{\eqref{y-sep}}{=}p_{k;t}^{b_{ki;t}}y_{i;t}y_{k;t}^{[b_{ki;t}]_+}(1\op y_{k;t})^{-b_{ki;t}}.
\end{align*}
Thus, the formula \eqref{y-sep} is proved. Since the mutation of the $\hy$-variables is the same as that in normalized cluster patterns, the separation formula \eqref{hy-sep} is also the same (\cite[Proposition 3.13]{CA4}). Finally, we show the formula \eqref{x-sep}. For $i\neq k$, the $G$-matrices, $F$-polynomials, and the $J$-functions are invariant under mutation; thus \eqref{x-sep} holds. For $i=k$, we have
\begin{align*}
x_{k;t'}&=\left(\prod_{j=1}^nx_j^{g_{jk;t'}}\right)\dfrac{F_{k;t'}(\hby)}{J_{k;t'}|_\PP(\by,\CP)}\\
&=\prod_{j=1}^n x_j^{-g_{jk;t}+\sum_{\ell=1}^n g_{j\ell;t}[-b_{\ell k;t}]_+-\sum_{\ell=1}^n b_{j\ell;t_0}[-c_{\ell k;t}]_+}\\
&\quad\times \dfrac{F_{k;t}(\hby)^{-1}\left(\displaystyle\prod_{j=1}^n\hy_j^{[-c_{jk;t}]_+}\right)\displaystyle\prod_{j=1}^nF_{j;t}(\hby)^{[-b_{jk;t}]_+}}{p_{k;t}^{-1}J_{k;t}|_\PP(\by,\CP)^{-1}\left(\displaystyle\prod_{j=1}^ny_j^{[-c_{jk;t}]_+}\right)\displaystyle\prod_{j=1}^nJ_{j;t}|_\PP(\by,\CP)^{[-b_{jk;t}]_+}}\\
&\quad\times \dfrac{1+\left(\displaystyle\prod_{j=1}^n\hy_j^{c_{jk;t}}\right)\displaystyle\prod_{j=1}^nF_{j;t}(\hby)^{b_{jk;t}}}{1\op \left(\displaystyle\prod_{j=1}^ny_j^{c_{jk;t}}\right)\displaystyle\prod_{j=1}^nJ_{j;t}|_\PP(\by,\CP)^{b_{jk;t}}}\\
&=p_{k;t}\prod_{j=1}^n x_j^{-g_{jk;t}}\dfrac{F_{k;t}(\hby)^{-1}}{J_{k;t}|_\PP(\by,\CP)^{-1}}\times\prod_{j=1}^n x_j^{\sum_{\ell=1}^n g_{j\ell;t}[-b_{\ell k;t}]_+}\dfrac{\displaystyle\prod_{j=1}^nF_{j;t}(\hby)^{[-b_{jk;t}]_+}}{\displaystyle\prod_{j=1}^nJ_{j;t}|_\PP(\by,\CP)^{[-b_{jk;t}]_+}}\\
&\quad\times\left(\prod_{j=1}^n x_j^{-\sum_{\ell=1}^n b_{j\ell;t_0}[-c_{\ell k;t}]_+}\right)\times \prod_{j=1}^n\prod_{\ell=1}^nx_{\ell}^{b_{\ell j;t_0}[-c_{jk;t}]_+}\\
&\quad\times\dfrac{1+\left(\displaystyle\prod_{j=1}^n\hy_j^{c_{jk;t}}\right)\displaystyle\prod_{j=1}^nF_{j;t}(\hby)^{b_{jk;t}}}{1\op \left(\displaystyle\prod_{j=1}^ny_j^{c_{jk;t}}\right)\displaystyle\prod_{j=1}^nJ_{j;t}|_\PP(\by,\CP)^{b_{jk;t}}}\\
&\overset{\eqref{x-sep}}{=}p_{k;t}x_{k;t}^{-1}\left(\prod_{j=1}^nx_{j;t}^{[-b_{jk;t}]_+}\right)\frac{1+\hy_{k;t}}{1\op y_{k;t}},
\end{align*}
where in the last equality, we used the result of \eqref{y-sep}, \eqref{hy-sep}.
\end{proof}
\subsection{Example of $J$-functions}
Let us present the $J$-functions of type $A_2$. Let $B_0=\begin{pmatrix}
0 & -1\\
1 & 0
\end{pmatrix}$, and we consider the case where the $J$-functions mutate in directions 1 and 2 alternately. 
\begin{align*}
&\begin{cases}
J_{1;0}(\by,\CP)=1\\
J_{2;0}(\by,\CP)=1
\end{cases},\quad
\begin{cases}
J_{1;1}(\by,\CP)=p_{1;0}^{-1}(1+y_1)\\
J_{2;1}(\by,\CP)=1
\end{cases},\quad
\begin{cases}
J_{1;2}(\by,\CP)=p_{1;0}^{-1}(1+y_1)\\
J_{2;2}(\by,\CP)=p_{2;1}^{-1}(1+p_{1;0}^{-1}y_2(1+y_1))
\end{cases},\\
&\begin{cases}
J_{1;3}(\by,\CP)=p_{1;2}^{-1}p_{1;0}y_1\dfrac{1+ p_{2;1}^{-1}y_1^{-1}+ p_{1;0}^{-1}p_{2;1}^{-1}y_1^{-1}y_2+ p_{1;0}^{-1}p_{2;1}^{-1}y_2}{1+ y_1}\\
J_{2;3}(\by,\CP)=p_{2;1}^{-1}(1+p_{1;0}^{-1}y_2(1+y_1))
\end{cases},\\
&\begin{cases}
J_{1;4}(\by,\CP)=p_{1;2}^{-1}p_{1;0}y_1\dfrac{1+ p_{2;1}^{-1}y_1^{-1}+ p_{1;0}^{-1}p_{2;1}^{-1}y_1^{-1}y_2+ p_{1;0}^{-1}p_{2;1}^{-1}y_2}{1+ y_1}\\
{}\\
J_{2;4}(\by,\CP)=p_{2;3}^{-1}p_{2;1}y_1y_2\dfrac{1+p_{1;2}^{-1}p_{1;0}y_2^{-1}\dfrac{1+ p_{2;1}^{-1}y_1^{-1}+ p_{1;0}^{-1}p_{2;1}^{-1}y_1^{-1}y_2+ p_{1;0}^{-1}p_{2;1}^{-1}y_2}{1+ y_1}}{1+p_{1;0}^{-1}y_2(1+y_1)}
\end{cases},\\
&\begin{cases}
J_{1;5}(\by,\CP)=p_{1;4}^{-1}p_{1;2}p_{1;0}^{-1}y_1^{-1}y_2\dfrac{1+ y_1}{1+ p_{2;1}^{-1}y_1^{-1}+ p_{1;0}^{-1}p_{2;1}^{-1}y_1^{-1}y_2+ p_{1;0}^{-1}p_{2;1}^{-1}y_2}\\[2ex]
\phantom{J_{1;5}(\by,\CP)=}\times \left(1+p_{2;3}^{-1}p_{2;1}y_1\dfrac{1+p_{1;2}^{-1}p_{1;0}y_2^{-1}\dfrac{1+ p_{2;1}^{-1}y_1^{-1}+ p_{1;0}^{-1}p_{2;1}^{-1}y_1^{-1}y_2+ p_{1;0}^{-1}p_{2;1}^{-1}y_2}{1+ y_1}}{1+p_{1;0}^{-1}y_2(1+y_1)}\right)\\[3ex]
J_{2;5}(\by,\CP)=p_{2;3}^{-1}p_{2;1}y_1y_2\dfrac{1+p_{1;2}^{-1}p_{1;0}y_2^{-1}\dfrac{1+ p_{2;1}^{-1}y_1^{-1}+ p_{1;0}^{-1}p_{2;1}^{-1}y_1^{-1}y_2+ p_{1;0}^{-1}p_{2;1}^{-1}y_2}{1+ y_1}}{1+p_{1;0}^{-1}y_2(1+y_1)}
\end{cases}.
\end{align*}
As we see, the $J$-functions are not polynomials in general.

\section{\texorpdfstring{$\CP$}{P}-normalized generalized cluster patterns}\label{P-GCP}

In this section, we define generalized cluster patterns following Chekhov and Shapiro \cite{CS14}, and introduce $\CP$-normalized generalized cluster patterns. We then observe that the mutation formulas for the $x$-, $y$-, and $\hy$-variables are parallel to those for ordinary cluster patterns.

\subsection{Generalized cluster patterns}

In the following, we provide the definition and the basic properties of generalized cluster patterns introduced by Chekhov and Shapiro \cite{CS14}. We follow the convention of exchange matrix $B$ by \cite{CA4} which is the transpose of \cite{CS14}.

Let $\br=(r_1,\dots,r_n)$ be an $n$-tuple of positive integers. We call each $r_i$ a \emph{polynomial degree}. Also, let $R=\diag(r_1,\dots,r_n)$ be a diagonal matrix with $r_i$ at the $(i,i)$-entry.

As in Section \ref{reform}, let $\PP$ be a semifield and $\CF$ be an ambient field.

\begin{defi} [Generalized seeds]
A \emph{(generalized) seed} is a triplet $\Sigma^p=(\bx,\bp,B)$ given by the following.
\begin{enumerate}
\item $\bx=(x_1,\dots,x_n)$ is an $n$-tuple of algebraically independent elements over $\Q\PP$ generating $\CF$, and each $x_i$ is called a \emph{cluster variable} or an \emph{$x$-variable}.
\item $\bp=(\bp_1,\dots,\bp_n)$ is an $n$-tuple where each $\bp_i=(p_{i0},\dots,p_{ir_i})$ is an $(r_i+1)$-tuple of elements in $\PP$. Each $p_{i\ell}$ is called a \emph{$p$-variable}.
\item $B$ is a skew-symmetrizable integer matrix with a skew-symmetrizer $D$, called an \emph{exchange matrix.}
\end{enumerate}
\end{defi}

\begin{defi} [Seed mutations] A seed $(\Sigma^p)'=(\bx',\bp',B')$ is called a \emph{mutation} of $\Sigma^p$ in direction $k=1,\dots,n$ if it satisfies the following rules:
\begin{align}
x_i'&=\begin{cases}
x_k^{-1}\left(\prod_{j=1}^nx_j^{[-b_{jk}]_+}\right)^{r_k}\sum_{\ell=0}^{r_k}p_{k\ell}\left(\prod_{j=1}^nx_j^{b_{jk}}\right)^\ell & \text{$(i=k)$},\\
x_i &\text{$(i\neq k)$},
\end{cases}\label{GCA-x-mut}\\
p_{k\ell}'&=p_{k,r_k-\ell},\label{GCA-p-mut1}\\
p_{i\ell}'/p_{i0}'&=\begin{cases}
(p_{kr_k})^{\ell b_{ki}}p_{i\ell}/p_{i0} & \text{$(i\neq k, b_{ki}\ge0)$},\\
(p_{k0})^{\ell b_{ki}}p_{i\ell}/p_{i0} & \text{$(i\neq k, b_{ki}<0)$},
\end{cases}\label{GCA-p-mut2}\\
b_{ij}'&=\begin{cases}
-b_{ij} & \text{$(i=k$ or $j=k)$},\\
b_{ij}+r_kb_{ik}[b_{kj}]_++r_k[-b_{ik}]_+b_{kj} & \text{$(i,j\neq k)$}.
\end{cases}\label{GCA-B-mut}
\end{align}
Now, let $\hat{B}=BR=(\hat{b}_{ij})$. Then, $\hat{B}$ mutates as follows:
\begin{align}
\hat{b}_{ij}'=\begin{cases}
-\hat{b}_{ij} & \text{$(i=k$ or $j=k)$},\\
\hat{b}_{ij}+\hat{b}_{ik}[\hat{b}_{kj}]_++[-\hat{b}_{ik}]_+\hat{b}_{kj} & \text{(otherwise)}.
\end{cases}\label{GCA-hB-mut}
\end{align}
From \eqref{GCA-hB-mut}, the mutation of $\hat{B}$ coincides with the mutation of the exchange matrix in cluster patterns, and $\hat{B}$ is a skew-symmetrizable matrix with a skew-symmetrizer $R D$.
\end{defi}

As in the ordinary case, $\bp'$ is not uniquely determined from $\bp$ by \eqref{GCA-p-mut1} and \eqref{GCA-p-mut2}.

\begin{defi} [Generalized cluster patterns] We attach a seed $\Sigma^p_t=(\bx_t,\bp_t,B_t)$ to each vertex $t \in \mathbb{T}_n$. The family of the seeds $\bsig^p=(\Sigma^p_t)_{t\in\T_n}$ is called a \emph{(non-normalized) generalized cluster pattern} if, for any two vertices $t$ and $t'$ connected by an edge labeled $k$, the seed $\Sigma^p_{t'}$ is the mutation of $\Sigma^p_{t}$ in direction $k$.
\end{defi}

The following is an example of the fact that many important properties of ordinary cluster patterns are preserved in generalized cluster patterns.

\begin{thm}[{\cite[Theorem 2.5]{CS14}}] \label{GCA-Laurent-pheno} Any cluster variable $x_{i;t}$ is expressed as a Laurent polynomial in the initial cluster variables with coefficients in $\Z\PP$.
\end{thm}

\subsection{\texorpdfstring{$\CP$}{P}-normalized generalized cluster patterns}
As in the ordinary case, a non-normalized generalized cluster pattern is not uniquely determined by an initial seed, in general. In this subsection, we introduce the $\CP$-normalization condition. Under the condition, a generalized cluster pattern is uniquely determined by an initial seed as in the ordinary case.

\begin{defi}[$\CP$-normalization condition]\label{GCA_P-norm cond}   Let $\CP$ be a $P$-pattern as in Definition \ref{P-pattern}. Suppose that a generalized cluster pattern $\bsig^p=(\Sigma_t^p)_{t\in\T_n}$ satisfies the following condition, which is called the \emph{$\CP$-normalization condition}:
\begin{align}
\bigoplus_{\ell=0}^{r_i}p_{i\ell;t}=p_{i;t}\label{GCAP-nor}
\end{align}
for any $i\in\{1,\dots,n\}$ and $t\in\T_n$. Such a pattern $\bsig^p$ is called a \emph{$\CP$-normalized generalized cluster pattern}. A $\CP$-normalized generalized cluster pattern is said to satisfy the \emph{normalization condition} if $\CP$ is trivial. Also, a $\CP$-normalized generalized cluster  pattern satisfying this condition is called a \emph{normalized generalized cluster pattern}.
\end{defi}

Let $\bsig^p$ be a $\CP$-normalized generalized cluster pattern. We introduce $y$-variables in parallel with ordinary case. For $i=1,\dots,n$ and $\ell=0,\dots,r_i$, we define
\begin{align}
y_{i\ell;t}=\dfrac{p_{i\ell;t}}{p_{i0;t}} \label{GCAydef}.
\end{align}

By the definition, we have $p_{i\ell;t}=p_{i0;t}y_{i\ell;t}$ for any $i=1,\dots,n$, $\ell=1,\dots,r_i$, and $t\in\T_n$. Combining it with \eqref{GCAP-nor}, we obtain $p_{i0;t}(1\op y_{i1;t}\op\dots\op y_{ir_i;t})=p_{i;t}$; therefore, we have
\begin{align}
p_{i0;t}=\dfrac{p_{i;t}}{1\op y_{i1;t}\op\cdots\op y_{ir_i;t}}.\label{GCA:p to y1}
\end{align}
Moreover, by \eqref{GCAydef} and \eqref{GCA:p to y1}, we obtain
\begin{align}
p_{i\ell;t}=\dfrac{p_{i;t}y_{i\ell;t}}{1\op y_{i1;t}\op\cdots\op y_{ir_i;t}}. \label{GCA:p to y2}
\end{align}

For each $i\in\{1,\dots,n\}$, let $\by_{i;t}=(y_{i1;t},\dots,y_{ir_i;t})$, and write $\by_t=(\by_{1;t},\dots,\by_{n;t})$. Each $y_{i\ell;t}$ is called a \emph{$y$-variable}.

\begin{prop}\label{prop:GCA-y-mut} Let $\bsig^p=(\Sigma^p_t)_{t\in\T_n}$ be a $\CP$-normalized generalized cluster pattern, and let $t'$ be $k$-adjacent to $t$. Then, we obtain the following.

\textup{(1)} The $y$-variable $y_{i\ell;t}$ mutates in direction $k$ as follows:
\begin{align}
y_{i\ell;t'}&=\begin{cases}
\dfrac{y_{k,r_k-\ell;t}}{y_{kr_k;t}} & \text{$(i=k)$},\\
p_{k;t}^{\ell b_{ki;t}}y_{i\ell;t}y_{kr_k;t}^{[\ell b_{ki;t}]_+}(1\op y_{k1;t}\op\cdots\op y_{kr_k;t})^{-\ell b_{ki;t}} & \text{$(i\neq k)$}.
\end{cases}\label{GCA-y-mut}
\end{align}

\textup{(2)} The mutation of a seed $\Sigma^p_t$ in direction $k$ is uniquely determined.
\end{prop}

\begin{proof}
The first case of \eqref{GCA-y-mut} immediately follows from \eqref{GCA-p-mut1}. We consider the second case of \eqref{GCA-y-mut}. First, we assume $b_{ki;t}\ge0$. Then, using the first case of \eqref{GCA-p-mut2}, \eqref{GCA:p to y1}, and \eqref{GCA:p to y2}, we calculate $y_{i\ell;t'}$ for $i\neq k$ as follows:
\begin{align*}
y_{i\ell;t'}&=y_{i\ell;t}\left(\dfrac{p_{k;t}y_{kr_k;t}}{1\op y_{k1;t}\op\cdots\op y_{kr_k;t}}\right)^{\ell b_{ki;t}}\\
&=p_{k;t}^{\ell b_{ki;t}}y_{i\ell;t}y_{kr_k;t}^{\ell b_{ki;t}}(1\op y_{k1;t}\op\cdots\op y_{kr_k;t})^{-\ell b_{ki;t}}.
\end{align*}
Since $[b_{ki;t}]_+=b_{ki;t}$, the claim holds for $i \neq k$ and $b_{ki;t}\ge0$. Next, we suppose $b_{ki;t}<0$. Then, we calculate as follows:
\begin{align*}
y_{i\ell;t'}&=y_{i\ell;t}\left(\dfrac{p_{k;t}}{1\op y_{k1;t}\op\cdots\op y_{kr_k;t}}\right)^{\ell b_{ki;t}}\\
&=p_{k;t}^{\ell b_{ki;t}}y_{i\ell;t}(1\op y_{k1;t}\op\cdots\op y_{kr_k;t})^{-\ell b_{ki;t}}.
\end{align*}
Since $[b_{ki;t}]_+=0$, the proof of (1) is complete.

(2) As in the case of ordinary cluster patterns, the claim holds by (1), \eqref{GCA:p to y1}, and \eqref{GCA:p to y2}.
\end{proof}

\begin{prop} A $\CP$-normalized generalized cluster pattern is uniquely determined by the $P$-pattern and an initial seed $\Sigma_{t_0}^p$. 
\end{prop}
\begin{proof}
It immediately follows from Proposition \ref{prop:GCA-y-mut} (2).
\end{proof}

Next, for $i=1,\dots,n$ and $\ell=1,\dots,r_i$, we define the \emph{$\hy$-variables} as follows:
\begin{align}
\hy_{i\ell;t}=y_{i\ell;t}\left(\prod_{j=1}^nx_{j;t}^{b_{j i;t}}\right)^\ell. \label{GCA-hy-def}
\end{align}

Then, \eqref{GCA-x-mut} is expressed using $y$- and $\hy$-variables as follows:

\begin{align}
x_{i;t'}=\begin{cases}
p_{k;t}x_{k;t}^{-1}\displaystyle\left(\prod_{j=1}^nx_{j;t}^{[-b_{jk;t}]_+}\right)^{r_k}\dfrac{1+\hy_{k1;t}+\dots+\hy_{kr_k;t}}{1\op y_{k1;t}\op\dots\op y_{kr_k;t}} & \text{$(i=k)$},\\
x_{i;t} & \text{$(i\neq k)$}.
\end{cases}\label{GCA-x-mut2}
\end{align}
Also, noting that $\hat{b}_{jk;t}=r_kb_{jk;t}$, the first case of \eqref{GCA-x-mut2} can also be expressed using the entries of $\hat{B}_t$:
\begin{align}
p_{k;t}x_{k;t}^{-1}\displaystyle\left(\prod_{j=1}^nx_{j;t}^{[-\hat{b}_{jk;t}]_+}\right)\dfrac{1+\hy_{k1;t}+\dots+\hy_{kr_k;t}}{1\op y_{k1;t}\op\dots\op y_{kr_k;t}}.
\end{align}

\begin{prop} The $\hy$-variables mutate in direction $k$ as follows:
\begin{align}
\hy_{i\ell;t'}&=\begin{cases}
\dfrac{\hy_{k,r_k-\ell;t}}{\hy_{kr_k;t}} & \text{$(i=k)$},\vspace{0.5em}\\
\hy_{i\ell;t}\hy_{kr_k;t}^{[\ell b_{ki;t}]_+}(1+\hy_{k1;t}+\cdots+\hy_{kr_k;t})^{-\ell b_{ki;t}} & \text{$(i\neq k)$}.
\end{cases}\label{GCA-hy-mut}
\end{align}
\end{prop}
\begin{proof}
By the second case of \eqref{GCA-x-mut2}, we have
\begin{align*}
\hy_{i\ell;t'}&=y_{i\ell;t'}\left(\prod_{\substack{j=1\\j\neq k}}^nx_{j;t}^{b_{ji;t'}}\right)^\ell x_{k;t'}^{\ell b_{ki;t'}}.
\end{align*}
First, for $i=k$, the claim immediately holds. Next, we show the second case of \eqref{GCA-hy-mut}. For $i\neq k$, by the mutations of the $x$-, $y$-variables, and the exchange matrix \eqref{GCA-x-mut2}, \eqref{GCA-y-mut}, \eqref{GCA-B-mut}, we calculate as follows:
\begin{align*}
\hy_{i\ell;t'}&=p_{k;t}^{\ell b_{ki;t}}y_{i\ell;t}y_{kr_k;t}^{[\ell b_{ki;t}]_+}(1\op y_{k1;t}\op\cdots\op y_{kr_k;t})^{-\ell b_{ki;t}}\\
&\quad \times \prod_{\substack{j=1\\j\neq k}}^nx_{j;t}^{\ell(b_{ji;t}+r_kb_{jk;t}[b_{ki;t}]_++r_k[-b_{jk;t}]_+b_{ki;t})}\\
&\quad \times \left\{p_{k;t}x_{k;t}^{-1}\displaystyle\left(\prod_{j=1}^nx_{j;t}^{[-b_{jk;t}]_+}\right)^{r_k}\dfrac{1+\hy_{k1;t}+\dots+\hy_{kr_k;t}}{1\op y_{k1;t}\op\dots\op y_{kr_k;t}}\right\}^{-\ell b_{ki;t}}\\
&=y_{i\ell;t}\left(\prod_{j=1}^nx_{j;t}^{b_{ji;t}}\right)^\ell \times y_{kr_k;t}^{[\ell b_{ki;t}]_+}\left(\prod_{\substack{j=1\\j\neq k}}^nx_{j;t}^{b_{jk;t}}\right)^{r_k[\ell b_{ki;t}]_+}\times (1+\hy_{k1;t}+\cdots+\hy_{kr_k;t})^{-\ell b_{ki;t}}\\
&=\hy_{i\ell;t}\hy_{kr_k;t}^{[\ell b_{ki;t}]_+}(1+\hy_{k1;t}+\cdots+\hy_{kr_k;t})^{-\ell b_{ki;t}}.
\end{align*}
This completes the proof.
\end{proof}

Note that the mutation of $\hy$-variables does not depend on $\CP$. Also, similar to ordinary cluster patterns, when $\CP$ is the trivial $P$-pattern, the mutation of $\hy$-variables coincides with that of $y$-variables.

\subsection{\texorpdfstring{$\vep$}{epsilon}-expressions}
As in the ordinary cases \cite[Proposition 2.7]{Nak23}, the $\vep$-expressions of mutations present as follows.

\begin{prop} [$\vep$-expressions]
Let $\bsig$ be a $\CP$-normalized generalized cluster pattern. For any $i\in\{1,\dots,n\}$ and $\ell\in\{1,\dots,r_i\}$, let $y_{i\ell;t}^{(1)}=y_{i\ell;t}, y_{i\ell;t}^{(-1)}=\dfrac{y_{i,r_i-\ell;t}}{y_{ir_i;t}}$, $\hy_{i\ell;t}^{(1)}=\hy_{i\ell;t}, \hy_{i\ell;t}^{(-1)}=\dfrac{\hy_{i,r_i-\ell;t}}{\hy_{ir_i;t}}$. Then, for any sign $\vep\in\{1,-1\}$, the mutations in direction $k$ \eqref{GCA-x-mut2}, \eqref{GCA-y-mut}, and \eqref{GCA-hy-mut} are expressed as follows respectively:

\begin{align}
x_{i;t'}&=\begin{cases}
p_{k;t}x_{k;t}^{-1}\displaystyle\left(\prod_{j=1}^nx_{j;t}^{[-\vep b_{jk;t}]_+}\right)^{r_k}\dfrac{1+\hy_{k1;t}^{(\vep)}+\dots+\hy_{kr_k;t}^{(\vep)}}{1\op y_{k1;t}^{(\vep)}\op\dots\op y_{kr_k;t}^{(\vep)}} & \text{$(i=k)$},\\
x_{i;t} & \text{$(i\neq k)$},
\end{cases}\label{GCA-x-mut-ep}\\
y_{i\ell;t'}&=\begin{cases}
y_{i\ell;t}^{(-1)} & \text{$(i=k)$},\\
p_{k;t}^{\ell b_{ki;t}}y_{i\ell;t}y_{kr_k;t}^{[\vep\ell b_{ki;t}]_+}(1\op y_{k1;t}^{(\vep)}\op\cdots\op y_{kr_k;t}^{(\vep)})^{-\ell b_{ki;t}} & \text{$(i\neq k)$},
\end{cases}\label{GCA-y-mut-ep}\\
\hy_{i\ell;t'}&=\begin{cases}
\hy_{i\ell;t}^{(-1)} & \text{$(i=k)$},\\
\hy_{i\ell;t}\hy_{kr_k;t}^{[\vep\ell b_{ki;t}]_+}(1+ \hy_{k1;t}^{(\vep)}+\cdots+\hy_{kr_k;t}^{(\vep)})^{-\ell b_{ki;t}} & \text{$(i\neq k)$}.
\end{cases}\label{GCA-y-mut-ep2}
\end{align}
\end{prop}
\begin{proof}

By the definitions of $y_{k\ell;t}^{(-1)}$ and $\hy_{k\ell;t}^{(-1)}$, we have
\begin{align*}
1\op y_{k1;t}^{(-1)}\op\cdots\op y_{kr_k;t}^{(-1)}&=y_{kr_k;t}^{-1}(1\op y_{k1;t}\op\cdots\op y_{kr_k;t}),\\
1+\hy_{k1;t}^{(-1)}+\cdots+\hy_{kr_k;t}^{(-1)}&=\hy_{kr_k;t}^{-1}(1+\hy_{k1;t}+\cdots+\hy_{kr_k;t}).
\end{align*}
By dividing the expression for $\vep=1$ by that for $\vep=-1$, we obtain
\begin{align*}
&\quad\left(\prod_{j=1}^nx_{j;t}^{[-b_{jk;t}]_+}\right)^{r_k}\dfrac{1+\hy_{k1;t}^{(1)}+\dots+\hy_{kr_k;t}^{(1)}}{1\op y_{k1;t}^{(1)}\op\dots\op y_{kr_k;t}^{(1)}}\times\left\{\left(\prod_{j=1}^nx_{j;t}^{[b_{jk;t}]_+}\right)^{r_k}\dfrac{1+\hy_{k1;t}^{(-1)}+\dots+\hy_{kr_k;t}^{(-1)}}{1\op y_{k1;t}^{(-1)}\op\dots\op y_{kr_k;t}^{(-1)}}\right\}^{-1}\\
&=\left(\prod_{j=1}^nx_{j;t}^{-b_{jk;t}}\right)^{r_k}\dfrac{\hy_{kr_k;t}}{y_{kr_k;t}}=1,\\
&\quad\quad y_{kr_k;t}^{[\ell b_{ki;t}]_+}(1\op y_{k1;t}^{(1)}\op\cdots\op y_{kr_k;t}^{(1)})^{-\ell b_{ki;t}}\left\{y_{kr_k;t}^{[-\ell b_{ki;t}]_+}(1\op y_{k1;t}^{(-1)}\op\cdots\op y_{kr_k;t}^{(-1)})^{-\ell b_{ki;t}}\right\}^{-1}\\
&=y_{kr_k;t}^{[\ell b_{ki;t}]_+-[-\ell b_{ki;t}]_+-\ell b_{ki;t}}=1.
\end{align*}
The argument for the $\hy$-variables is similar to that for the $y$-variables; thus this completes the proof.
\end{proof}

\section{Separation formulas for generalized cluster patterns}\label{GCA-Sep}

In this section, we define the $C$-, $G$-matrices, and $F$-polynomials for generalized cluster patterns in parallel with those for ordinary cluster patterns. Then, using them, we provide the separation formulas for normalized and $\CP$-normalized generalized cluster patterns.

\subsection{\texorpdfstring{$C$}{C}-, \texorpdfstring{$G$}{G}-matrices, and \texorpdfstring{$F$}{F}-polynomials}

We introduce the generalized $C$-, $G$-matrices, and the generalized $F$-polynomials using the information of the generalized cluster patterns with principal coefficients, and observe their mutations. In the following, we consider normalized generalized cluster patterns.

\begin{defi}[Generalized cluster patterns with principal coefficients] Let $\Tp(\bu)$ be the set of all monomials in $\bu=(u_1,\dots,u_m)$, endowed with the standard multiplication and with the addition called the tropical sum given by
\begin{align}
\prod_{j=1}^mu_j^{a_j}\op\prod_{j=1}^mu_j^{b_j}=\prod_{j=1}^mu_j^{\min \{a_j,b_j\}}\label{Tp-sum}.
\end{align}
This forms a semifield, called a tropical semifield. For an initial vertex $t_0\in\T_n$, a generalized cluster pattern is called a \emph{generalized cluster pattern with principal coefficients} if its coefficient semifield is the tropical semifield $\Tp(\by)$, where $\bu$ is replaced by $\by=\by_{t_0}$. \vv
\end{defi}

\begin{defi}[$C$-matrices]\label{Cdef} Let $ r =\sum_{j=1}^nr_j$. For a generalized cluster pattern with principal coefficients, we define an $ r \times r $ matrix $C_t=(c_{m\ell;t}^{(ij)})$ called a \emph{(generalized) $C$-matrix} as follows: 
\begin{align}
y_{j\ell;t}&=\prod_{i=1}^n\prod_{m=1}^{r_i}y_{im}^{c_{m\ell;t}^{(ij)}},\label{GCA-C-def}
\end{align}
where $c_{m\ell;t}^{(ij)}$ is the $(m,\ell)$-entry of the $(i,j)$-block when $C_t$ is partitioned into $n \times n$ blocks of sizes $r_1, \dots, r_n$. Moreover, the $\ell$-th column vector $\bc_{\ell;t}^{(ij)}$ of the $(i,j)$-block matrix of a $C$-matrix $C_t$ is called a \emph{(generalized) $c$-vector}, and the family of $C$-matrices $\{C_t\}_{t\in\T_n}$ is called the \emph{generalized $C$-pattern}. We also set $c_{m0;t}^{(ij)}=0$ for convenience.
\end{defi}

Then, we have the mutation formula of $C$-matrices.

\begin{prop}\label{prop:GCA-C-mut} The $C$-matrices are uniquely determined by the following recurrence relations. The transformation $C_{t}\mapsto C_{t'}$ given by these relations is called the mutation in direction $k$ of $C$-matrices.
\begin{align}
C_{t_0}=I_{ r},\hh
c_{m\ell;t'}^{(ij)}=
\begin{cases}
c_{m,r_k-\ell;t}^{(ik)}-c_{mr_k;t}^{(ik)} & \text{$(j=k)$},\\
c_{m\ell;t}^{(ij)}+c_{mr_k;t}^{(ik)}[\ell b_{kj;t}]_++\left[\max_{1\le s\le r_k}\left\{-c_{ms;t}^{(ik)}\right\}\right]_+(\ell b_{kj;t})  & \text{$(j\neq k)$}.
\end{cases}\label{GCA-C-mut}
\end{align}
\end{prop}
\begin{proof}
The initial condition $C_{t_0}=I_r$ obviously holds. Let $t'$ be $k$-adjacent to $t$. By the definition of $C$-matrices \eqref{GCA-C-def} and the tropical sum \eqref{Tp-sum}, the mutation of $y$-variables in the tropical semifield $\Tp(\by)$ is given by
\begin{align}
y_{i\ell;t'}=\begin{cases}
\displaystyle\prod_{j=1}^n\prod_{m=1}^{r_j}y_{jm}^{c_{m,r_k-\ell;t}^{(jk)}}\left(\displaystyle\prod_{j=1}^n\prod_{m=1}^{r_j}y_{jm}^{c_{mr_k;t}^{(jk)}}\right)^{-1} & \text{$(i=k)$},\\[2em]
\displaystyle\prod_{j=1}^n\prod_{m=1}^{r_j}y_{jm}^{c_{m\ell;t}^{(ji)}} \prod_{j=1}^n\prod_{m=1}^{r_j}y_{jm}^{[\ell b_{ki;t}]_+c_{mr_k;t}^{(jk)}}\prod_{j=1}^n\prod_{m=1}^{r_j}y_{jm}^{(-\ell b_{ki;t})\min\{0,c_{m1;t}^{(jk)},\dots,c_{mr_k;t}^{(jk)}\}} & \text{$(i\neq k)$}.
\end{cases}
\end{align}
We note the following identity:
\begin{align}
-\min\{0,a_1,\dots,a_n\}=[\max\{-a_1,\dots,-a_n\}]_+.
\end{align}
By comparing the exponents of $y_{jm}$, we have \eqref{GCA-C-mut}.

\end{proof}

Next, we introduce $F$-polynomials. As a preparation, we show the following proposition, which is parallel to the ordinary case \cite[Proposition 11.2]{CA2}.

\begin{prop}\label{GCA-polynomialization} For any generalized cluster pattern with principal coefficients, $x_{i;t}\in\Z[\bx^{\pm1}, \by]$. 
\end{prop}
\begin{proof}
As in the proof of \cite[Proposition 11.2]{CA2}, we show the proposition. Below, we prove the following stronger statement.
\begin{claim}
For each $j$, $x_{i;t}$ is a polynomial in $y_{jm}$ whose constant term is a non-zero subtraction-free polynomial in $\bx^{\pm1}$ and $\by$ other than $y_{jm}$.
\end{claim}
We prove this claim by induction on $t\in\T_n$. For $t=t_0$, the claim clearly holds. Suppose that the claim holds for $t$, and let $t'\in\T_n$ be a vertex which is $k$-adjacent to $t$. The mutation of $x$-variables for generalized cluster patterns with principal coefficients is written using $C$-matrices as follows:
\begin{align}
x_{i;t'}=\begin{cases}
\displaystyle x_{k;t}^{-1}\left(\prod_{j=1}^nx_{j;t}^{[- b_{jk;t}]_+}\right)^{r_k}\sum_{\ell=0}^{r_k}\left\{\prod_{j=1}^n\prod_{m=1}^{r_j}y_{jm}^{[\max_{1\le s\le r_k}\{-c_{ms;t}^{(jk)}\}]_++c_{m\ell;t}^{(jk)}}\left(\prod_{j=1}^nx_{j;t}^{ b_{jk;t}}\right)^\ell\right\} & \text{$(i=k)$},\\
x_{i;t} & \text{$(i\neq k)$}.
\end{cases}\label{GCA-princ-x-mut}
\end{align}
First, for $i\neq k$, the claim holds by the induction hypothesis. Next, we consider the case $i=k$. For any pair $(j,m)$, at least one of the values
\begin{align}
\left[\max_{1\le s \le r_k}\{-c_{ms;t}^{(jk)}\}\right]_++c_{m\ell;t}^{(jk)} \quad(\ell=0,\dots,r_k) 
\end{align}
is equal to $0$.
Therefore, the numerator of \eqref{GCA-princ-x-mut} is a polynomial in $y_{jm}$ whose constant term is a polynomial in $\bx_{t_0}$ and $\by$ other than $y_{jm}$. 
The denominator of the first case of \eqref{GCA-princ-x-mut} $x_{k;t}$ also has the same properties by the induction hypothesis. By Theorem \ref{GCA-Laurent-pheno}, $x_{i;t'}\in\Z[\bx^{\pm1},\by^{\pm1}]$ holds. Thus, the numerator is divisible by the denominator, and the claim holds.
\end{proof}

\begin{defi}[$F$-polynomials] By Proposition \ref{GCA-polynomialization}, any cluster variable $x_{i;t}$ belongs to $\Z[\bx^{\pm1},\by]$. We define a (generalized) \emph{$F$-polynomial $F_{i;t}(\by)$} as the specialization of $x_{i;t}$, expressed as a Laurent polynomial in the initial $x$-variables $\bx$ and $y$-variables $\by$, to $x_1 = \dots = x_n = 1$. The family of $F$-polynomials $\{F_{i;t}(\by)\}_{t\in\T_n}$ is called the \emph{generalized $F$-pattern}.
\end{defi}

Then, the following recurrence relations hold.  

\begin{prop}\label{prop:GCA-F-mut} The $F$-polynomials are uniquely determined by the following recurrence relations. The transformation $F_{i;t}(\by)\mapsto F_{i;t'}(\by)$ given by these relations is called the mutation in direction $k$ of $F$-polynomials.
\begin{align}
F_{i;t_0}(\by)&=1,\\
F_{i;t'}(\by)&=\begin{cases}
\dfrac{M_{k;t}(\by)}{F_{k;t}(\by)} & \text{$(i=k)$},\\
F_{i;t}(\by) & \text{$(i\neq k)$},
\end{cases}\label{GCA-F-mut}
\end{align}
where
\begin{align}
\begin{split}
\displaystyle M_{k;t}(\by)&=\left(\prod_{j=1}^nF_{j;t}(\by)^{[-r_k b_{jk;t}]_+}\right)\left(\prod_{j=1}^n\prod_{m=1}^{r_j}y_{jm}^{\left[\max_{1\le s\le r_k}\left\{-c_{ms;t}^{(jk)}\right\}\right]_+}\right)\\
&\quad\times\left[1+
\sum_{\ell=1}^{r_k}\left\{\left(\prod_{j=1}^n\prod_{m=1}^{r_j}y_{jm}^{c_{m\ell;t}^{(jk)}}\right)\prod_{j=1}^nF_{j;t}(\by)^{\ell b_{jk;t}}\right\}\right].
\label{GCA-M-def}
\end{split}
\end{align}
\end{prop}
\begin{proof}
The initial condition $F_{i;t_0}(\by)=1$ clearly holds. Furthermore, replacing each $x_{i;t}$ with $F_{i;t}(\by)$ in the mutation of $x$-variables \eqref{GCA-princ-x-mut}, we obtain the claim.
\end{proof}

Let $\Q_{\mathrm{sf}}(\by)$ be the set of elements in the rational function field $\Q(\by)$ that have subtraction-free expressions. The semifield homomorphism $\pi_{\mathrm{trop}} : \Q_{\mathrm{sf}}(\by)\rightarrow \mathrm{Trop}(\by)$, $y_{i\ell}\mapsto y_{i\ell}$ is called the \emph{tropicalization homomorphism}. Since the recurrence relation for $F$-polynomials is subtraction-free, each $F_{i;t}(\by)$ is an element of $\Q_{\mathrm{sf}}(\by)$. Therefore, the homomorphism can be applied to $F_{i;t}(\by)$.

\begin{prop} $\pi_{\textup{trop}}(F_{i;t}(\by))=1$.
\end{prop}
\begin{proof}
We prove this by induction on $t\in\T_n$. The case $t=t_0$ obviously holds. Assuming that it holds for $t$, let $t'$ be the $k$-adjacent vertex to $t$. It suffices to show that $\pi_{\textup{trop}}(M_{k;t}(\by))=1$. By the assumption, in \eqref{GCA-M-def}, 
\begin{align*}
\pi_{\textup{trop}}\left(1+\sum_{\ell=1}^{r_k}\left\{\left(\prod_{j=1}^n\prod_{m=1}^{r_j}y_{jm}^{c_{m\ell;t}^{(jk)}}\right)\prod_{j=1}^nF_{j;t}(\by)^{\ell b_{jk;t}}\right\}\right)=\prod_{j=1}^n\prod_{m=1}^{r_j}y_{jm}^{-\left[\max_{1\le s\le r_k}\left\{-c_{ms;t}^{(jk)}\right\}\right]_+}.
\end{align*}
Thus, we have the claim.
\end{proof}

In the following, we introduce $G$-matrices. As a preparation, we define the principal $\Z^n$-grading.

\begin{defi} [Principal $\Z^n$-grading]
We define the \emph{principal} $\Z^n$\emph{-grading} $\deg$ on $\Z[\bx^{\pm1},\by]$ in parallel to \cite{CA4} as follows: 
\begin{align}
\deg(x_i)=\be_i,\hh \deg(y_{i\ell})=-\ell \bb_{i;t_0},\label{eq:deg}
\end{align}
where $\be_i$ is the $i$-th unit vector and $\bb_{i;t_0}$ is the $i$-th column vector of the initial exchange matrix $B_{t_0}$. 
\end{defi}

By Theorem \ref{GCA-Laurent-pheno}, any cluster variable $x_{i;t}$ is expressed as a Laurent polynomial in the initial cluster variables. By the definition of $\hy$-variables, $\deg(\hy_{i\ell;t_0})=\bze \hh(i=1,\dots,n,\ell=1,\dots,r_i)$. Also, as in the case of ordinary cluster patterns \cite[proof of Proposition 6.1]{CA4}, any cluster variable is homogeneous with respect to the principal $\Z^n$-grading. Therefore, we regard $\deg(\hy_{i\ell;t})=\bze$ for any $i\in\{1,\dots,n\}, \ell\in\{1,\dots,r_i\}, t\in\T_n$. We define the $G$-matrix as in the ordinary cluster patterns.

\begin{defi} [Generalized $G$-matrices] For a generalized cluster pattern with principal coefficients, we define a \emph{(generalized) $g$-vector} as follows:
\begin{align}
\bg_{i;t}=\deg(x_{i;t})\in\Z^n.
\end{align}
Let $G_t$ be the $n \times n$ integer matrix whose $i$-th column vector is $\bg_{i;t}$. This is called a \emph{(generalized) $G$-matrix} and the family of $G$-matrices $\{G_t\}_{t\in\T_n}$ is called the \emph{generalized $G$-pattern}.
\end{defi}

Then, we obtain the following mutation formula for $G$-matrices.

\begin{prop} The $G$-matrices are uniquely determined by the following recurrence relations. The transformation $G_{t}\mapsto G_{t'}$ given by these relations is called the mutation of $G$-matrices.
\begin{align}
G_{t_0}&=I_n,\\
g_{ij;t'}&=\begin{cases}
-g_{ik;t}+\displaystyle r_k\sum_{\ell=1}^ng_{i\ell;t}[-b_{\ell k;t}]_+-\displaystyle\sum_{\ell=1}^n\sum_{m=1}^{r_\ell}m b_{i\ell;t_0}\left[\max_{1\le s \le r_k}\left\{-c_{ms;t}^{(\ell k)}\right\}\right]_+ & \text{$(j=k)$},\\
g_{ij;t} & \text{$(j\neq k)$}.
\end{cases}\label{GCA-G-mut}
\end{align}
\end{prop}
\begin{proof}
The initial condition follows from \eqref{eq:deg}. By the definitions of the tropical sum and the principal $\Z^n$-grading, we have
\begin{align}
\deg \left(\dfrac{1}{1\op y_{k1;t} \op \cdots \op y_{kr_k;t}}\right)&=\deg \left(\prod_{\ell=1}^n\prod_{m=1}^{r_\ell}y_{\ell m}^{-\min\{0,c_{m1;t}^{(\ell k)},\dots,c_{mr_k;t}^{(\ell k)}\}}\right)\notag\\
&=\deg \left(\prod_{\ell=1}^n\prod_{m=1}^{r_\ell}y_{\ell m}^{[\max_{1\le s \le r_k}\{-c_{ms;t}^{(\ell k)}\}]_+}\right)\\
&=-\sum_{\ell=1}^n\sum_{m=1}^{r_\ell}m\bb_{\ell;t_0}\left[\max_{1\le s \le r_k}\{-c_{ms;t}^{(\ell k)}\}\right]_+.\label{deg-y-sum}
\end{align}

Here, we formally set $\deg(\hy_{i\ell;t})=0$, which is a rational function of initial $\hy$-variables $\hby$. Therefore, we have the claim by \eqref{deg-y-sum} and the mutation of $x$-variables \eqref{GCA-x-mut2}.
\end{proof}

\subsection{Block structure of $C$-matrices}

We investigate the structure of \texorpdfstring{$C$}{C}-matrices. In particular, we examine the structure of each block matrix.

The following notion is necessary to classify the block matrices of $C$-matrices.

\begin{defi} [Mutation signs] For an initial vertex $t_0\in\T_n$ and $i\in\{1,\dots,n\}, t\in\T_n$, we define $\sigma_{i;t}\in\{\pm1\}$ as follows: 
\begin{align}
\sigma_{i;t}=\begin{cases}
1 & \text{(the label $i$ appears an even number of times in $t$)},\\
-1 & \text{(the label $i$ appears an odd number of times in $t$)}.
\end{cases}
\end{align}
The sign $\sigma_{i;t}$ is called the \emph{mutation sign} at $(i,t)$. Note that the mutation sign $\sigma_{i;t}$ also depends on $t_0$.
\end{defi}

\begin{prop} \label{prop:non-nega} \textup{(1)} In each block matrix of a $C$-matrix, the entries except for those in the last row are nonnegative, and the entries in the last column except for the one in the last row are all zero.

\textup{(2)} Suppose that $m\neq r_i$. Then, the following holds.
\begin{align}
c_{m\ell;t'}^{(ij)}=\begin{cases}
c_{m,r_k-\ell;t}^{(ik)} & \text{$(j=k)$},\\
c_{m\ell;t}^{(ij)} & \text{$(j\neq k)$}.
\end{cases}
\end{align}
\end{prop}

\begin{proof}
We show the claim (1). We prove these claims simultaneously by the induction on $t\in\T_n$. For $t=t_0$, they hold because $C_{t_0}$ is the identity matrix. Assume that they hold for $t\in\T_n$, and let $t'\in\T_n$ be $k$-adjacent to $t$. First, for the $(i,k)$-block, by \eqref{GCA-C-mut}, we have
\begin{align*}
c_{m\ell;t'}^{(ik)} = c_{m,r_k-\ell;t}^{(ik)} - c_{mr_k;t}^{(ik)} = c_{m,r_k-\ell;t}^{(ik)} \ge 0 \quad (m\neq r_i).
\end{align*}
In particular, when $\ell=r_k$, we have $c_{m\ell;t'}^{(ik)}=0$. Next, we consider the entries of the $(i,j)$-block for $j\neq k$. For $m\neq r_i$, we obtain
\begin{align*}
c_{m\ell;t'}^{(ij)} &\overset{\eqref{GCA-C-mut}}{=} c_{m\ell;t}^{(ij)} + c_{mr_k;t}^{(ik)}[\ell b_{kj;t}]_+ + \left[\max_{1\le s \le r_k}\{-c_{ms;t}^{(ik)}\}\right]_+\ell b_{kj;t}\\
&= c_{m\ell;t}^{(ij)} \ge 0.
\end{align*}
In particular, when $\ell=r_j$, we have $c_{m\ell;t'}^{(ij)}=0$. Thus, we have the claim (1). The claim (2) is proved by the proof of (1).
\end{proof}

We define the $n \times n$ square matrix $\tilde{I}_n$ obtained by reversing the order of the rows of the $n \times n$ identity matrix. That is, the $(i,j)$-entry $a_{ij}$ of $\tilde{I}_n$ is defined as follows:
\begin{align*}
a_{ij}=\begin{cases}
1 & \text{$(i+j=n+1)$},\\
0 & \text{(otherwise)}.
\end{cases}
\end{align*}
Then, we organize the block matrices of the $C$-matrices as follows.
\begin{prop}\label{Block-str}
The $(i,j)$-block matrix $C_{ij;t}$ of the $C$-matrix $C_t$ is expressed in the following form.
\begin{align*}
C_{ii;t}=\begin{cases}
\begin{pmatrix}
\multicolumn{3}{c}{I_{r_i-1}} & \bze \\
 \ast & \cdots & \ast & \ast
\end{pmatrix} & \text{$(\sigma_{i;t}=1)$},\vspace{0.5em}\\ \begin{pmatrix}
\multicolumn{3}{c}{\tilde{I}_{r_i-1}} & \bze \\
 \ast & \cdots & \ast & \ast
\end{pmatrix} & \text{$(\sigma_{i;t}=-1)$},
\end{cases}\quad
C_{ij;t}=\begin{pmatrix}
\multicolumn{3}{c}{O_{(r_{i}-1)\times r_j}}\\
\ast & \cdots & \ast 
\end{pmatrix} \quad \text{$(i\neq j)$}.
\end{align*}
\end{prop}
\begin{proof}
The claim is shown by Proposition \ref{prop:non-nega} (2).
\end{proof}

In the following, we consider the last row of the block matrices of $C$-matrices. As in the following proposition, the last row is uniquely determined by its last column entry.
\begin{prop}\label{last row of C} \textup{(1)} For $\ell=1,\dots,r_j-1$, the following holds. 
\begin{align}
c_{r_i\ell;t}^{(ij)}&=\begin{cases}
\dfrac{\ell}{r_j}c_{r_i r_j;t}^{(ij)} & \text{$(i\neq j)$},\\[1em]
\dfrac{\ell}{r_j}(c_{r_jr_j;t}^{(jj)}-1) & \text{$(i=j,\sigma_{j;t}=1)$},\\[1em]
\dfrac{\ell}{r_j}\left(c_{r_jr_j;t}^{(jj)}-\dfrac{r_j-\ell}{\ell}\right) & \text{$(i=j,\sigma_{j;t}=-1)$}.
\end{cases}\label{eq:last row of C}
\end{align}
\textup{(2)} The entries in the $r_i$-th row of the $(i,j)$-block of the $C$-matrix are either all nonnegative or all nonpositive, and satisfy $|c_{r_i\ell;t}^{(ij)}|\le|c_{r_i,\ell+1;t}^{(ij)}|$ for $\ell=1,\dots,r_j-1$.
\end{prop}
\begin{proof}
First, we prove the claim (1) by the induction on $t$. It is obvious for $t=t_0$. Assume it holds for $t$ and let $t'$ be $k$-adjacent to $t$. For $j=k$, 
\begin{align}
c_{r_i\ell;t'}^{(ik)}&\overset{\eqref{GCA-C-mut}}{=}c_{r_i,r_k-\ell;t}^{(ik)}-c_{r_i r_k;t}^{(ik)}\notag\\
&\overset{\eqref{eq:last row of C}}{=}\begin{cases}
\dfrac{r_k-\ell}{r_k}c_{r_ir_k;t}^{(ik)}-c_{r_ir_k;t}^{(ik)}=-\dfrac{\ell}{r_k}c_{r_ir_k;t}^{(ik)}=\dfrac{\ell}{r_k}c_{r_ir_k;t'}^{(ik)} & \text{$(i\neq k)$},\\[1em]
\dfrac{r_k-\ell}{r_k}(c_{r_kr_k;t}^{(kk)}-1)-c_{r_kr_k;t}^{(kk)}=\dfrac{\ell}{r_k}(c_{r_kr_k;t'}^{(kk)}-\dfrac{r_k-\ell}{\ell}) & \text{$(i=k,\sigma_{k;t}=1)$},\\[1em]
\dfrac{r_k-\ell}{r_k}\left(c_{r_kr_k;t}^{(kk)}-\dfrac{r_k-(r_k-\ell)}{r_k-\ell}\right)-c_{r_kr_k;t}^{(kk)}=\dfrac{\ell}{r_k}(c_{r_kr_k;t'}^{(kk)}-1)  & \text{$(i=k,\sigma_{k;t}=-1)$}.
\end{cases}
\end{align}
Noting that the mutation sign changes, the claim holds. For $j\neq k$, we have
\begin{align}
c_{r_i\ell;t'}^{(ij)}&\overset{\eqref{GCA-C-mut}}{=}c_{r_i\ell;t}^{(ij)}+c_{r_ir_k;t}^{(ik)}[\ell b_{kj;t}]_++\left[\max_{1\le s\le r_k}\left\{-c_{r_is;t}^{(ik)}\right\}\right]_+(\ell b_{kj;t})\notag\\
&\overset{\eqref{eq:last row of C}}{=}
\begin{cases}
\begin{aligned}
&\dfrac{\ell}{r_j}c_{r_ir_j;t}^{(ij)}+\dfrac{\ell}{r_j}c_{r_ir_k;t}^{(ik)}[r_jb_{kj;t}]_+\\[0.5em]
&\quad+\dfrac{\ell}{r_j}\left[\max_{1\le s\le r_k}\left\{-c_{r_is;t}^{(ik)}\right\}\right]_+r_jb_{kj;t}
\end{aligned} & \text{$(i\neq k)$}\\[2.5em]
\begin{aligned}
&\dfrac{\ell}{r_j}(c_{r_jr_j;t}^{(jj)}-1)+\dfrac{\ell}{r_j}c_{r_jr_k;t}^{(jk)}[r_jb_{kj;t}]_+\\[0.5em]
&\quad+\dfrac{\ell}{r_j}\left[\max_{1\le s\le r_k}\left\{-c_{r_js;t}^{(jk)}\right\}\right]_+r_jb_{kj;t}
\end{aligned} & \text{$(i=j,\sigma_{j;t}=1)$}\\[2.5em]
\begin{aligned}
&\dfrac{\ell}{r_j}\left(c_{r_jr_j;t}^{(jj)}-\dfrac{r_j-\ell}{\ell}\right)+\dfrac{\ell}{r_j}c_{r_jr_k;t}^{(jk)}[r_jb_{kj;t}]_+\\[0.5em]
&\quad+\dfrac{\ell}{r_j}\left[\max_{1\le s\le r_k}\left\{-c_{r_js;t}^{(jk)}\right\}\right]_+r_jb_{kj;t}
\end{aligned} & \text{$(i=j,\sigma_{j;t}=-1)$}
\end{cases}\notag\\
&\overset{\eqref{GCA-C-mut}}{=}
\begin{cases}
\dfrac{\ell}{r_j}c_{r_ir_j;t'}^{(ij)} & \text{$(i\neq j)$},\\[1em]
\dfrac{\ell}{r_j}(c_{r_jr_j;t'}^{(jj)}-1) & \text{$(i=j,\sigma_{j;t}=1)$},\\[1em]
\dfrac{\ell}{r_j}\left(c_{r_jr_j;t'}^{(jj)}-\dfrac{r_j-\ell}{\ell}\right) & \text{$(i=j,\sigma_{j;t}=-1)$}.
\end{cases}
\end{align}
Next, we show the claim (2) by induction on $t$. In the following, we say that $a$ and $b$ have the same sign if $a$ and $b$ are either both nonnegative or both nonpositive. 
In the first and second case of \eqref{eq:last row of C}, for $\ell\in\{1,\dots,r_j-2\}$, it holds $|c_{r_i\ell;t}^{(ij)}|\le|c_{r_i,\ell+1;t}^{(ij)}|$ clearly. Furthermore, in the first case of \eqref{eq:last row of C}, it holds $|c_{r_i\ell;t}^{(ij)}|\le|c_{r_ir_j;t}^{(ij)}|$. Also, in the second case of \eqref{eq:last row of C}, since $c_{r_jr_j;t}^{(jj)}=1+\dfrac{r_j}{\ell}c_{r_j\ell;t}^{(jj)}$ for $\ell\neq r_j$, $c_{r_jr_j;t}^{(jj)}$ and $c_{r_j\ell;t}^{(jj)}$ have the same sign. Because the $C$-matrix is an integer matrix, $c_{r_jr_j;t}^{(jj)}=1+\dfrac{r_j}{\ell}c_{r_j\ell;t}^{(jj)}=\left(1+\dfrac{r_j-\ell}{\ell}c_{r_j\ell;t}^{(jj)}\right)+c_{r_j\ell;t}^{(jj)}$ is an integer; thus $\dfrac{r_j-\ell}{\ell}c_{r_j\ell;t}^{(jj)}$ is an integer. Therefore, $1+\dfrac{r_j-\ell}{\ell}c_{r_j\ell;t}^{(jj)}$ and $c_{r_j\ell;t}^{(jj)}$ have the same sign; thus we have $|c_{r_jr_j;t}^{(jj)}|=\left|\left(1+\dfrac{r_j-\ell}{\ell}c_{r_j\ell;t}^{(jj)}\right)+c_{r_j\ell;t}^{(jj)}\right|\ge|c_{r_j\ell;t}^{(jj)}|$. We consider the third case of equation \eqref{eq:last row of C}. By the equation, we have 
\begin{align}
\ell r_j^{-1}(c_{r_jr_j;t}^{(jj)}+1)=c_{r_j\ell;t}^{(jj)}+1 \quad (\ell=1,\dots,r_j-1). \label{eq:third case}
\end{align}
By \eqref{eq:third case}, $c_{r_jr_j;t}^{(jj)}+1$ and $c_{r_j\ell;t}^{(jj)}+1$ have the same sign for $\ell=1,\dots,r_j$. Furthermore, we obtain that $c_{r_j\ell;t}^{(jj)}=-1$ is equivalent to $c_{r_jr_j;t}^{(jj)}=-1$, $c_{r_j\ell;t}^{(jj)}<-1$ is equivalent to $c_{r_jr_j;t}^{(jj)}<-1$, and $c_{r_j\ell;t}^{(jj)}\ge0$ is equivalent to $c_{r_jr_j;t}^{(jj)}\ge0$. Moreover, by \eqref{eq:third case}, we have $|c_{r_j\ell;t}^{(jj)}+1|\le |c_{r_j,\ell+1;t}^{(jj)}+1|$. By a similar argument, we obtain that $c_{r_j\ell;t}^{(jj)}+1$ and $c_{r_j,\ell+1;t}^{(jj)}+1$ have the same sign; thus we have $|c_{r_j\ell;t}^{(jj)}|\le |c_{r_j,\ell+1;t}^{(jj)}|$.
\end{proof}

\begin{prop}\label{prop:c-max} For any $i\in\{1,\dots,n\}, m\in\{1,\dots,r_i\},t\in\T_n$, the following holds:
\begin{align}
\left[\max_{1\le s\le r_k}\left\{-c_{ms;t}^{(ik)}\right\}\right]_+=[-c_{mr_k;t}^{(ik)}]_+. \label{cmax}
\end{align}
\end{prop}
\begin{proof}
It is immediately obtained from Proposition \ref{last row of C} (2).
\end{proof}

From the above result, the mutation formulas of the $C$-, $G$-matrices, and the $F$-polynomials \eqref{GCA-C-mut}, \eqref{GCA-G-mut}, and \eqref{GCA-M-def} are simplified as follows.

\begin{prop}\label{prop:new-mut} The recurrence relations for $C$-matrices, $G$-matrices, and $M_{k;t}(\by)$ appearing in the recurrence for $F$-polynomials are simplified as follows:
\begin{align}
c_{m\ell;t'}^{(ij)}&=
\begin{cases}
c_{m,r_k-\ell;t}^{(ik)}-c_{mr_k;t}^{(ik)} & \text{$(j=k)$},\\
c_{m\ell;t}^{(ij)}+c_{mr_k;t}^{(ik)}[\ell b_{kj;t}]_++[-c_{mr_k;t}^{(ik)}]_+(\ell b_{kj;t})  & \text{$(j\neq k)$},
\end{cases}\\
g_{ij;t'}&=\begin{cases}
-g_{ik;t}+\displaystyle\sum_{\ell=1}^ng_{i\ell;t}[-r_k b_{\ell k;t}]_+-\displaystyle\sum_{\ell=1}^n r_\ell b_{i\ell;t_0}\left[-c_{r_\ell r_k;t}^{(\ell k)}\right]_+ & \text{$(j=k)$},\\
g_{ij;t} & \text{$(j\neq k)$},
\end{cases}\label{GCA-G-mut2}\\
\begin{split}
\displaystyle M_{k;t}(\by)&=\left(\prod_{j=1}^nF_{j;t}(\by)^{[-r_k b_{jk;t}]_+}\right)\left(\prod_{j=1}^n\prod_{m=1}^{r_j}y_{jm}^{[-c_{mr_k;t}^{(jk)}]_+}\right) \\
&\quad\times\left[1+
\sum_{\ell=1}^{r_k}\left\{\left(\prod_{j=1}^n\prod_{m=1}^{r_j}y_{jm}^{c_{m\ell;t}^{(jk)}}\right)\prod_{j=1}^nF_{j;t}(\by)^{\ell b_{jk;t}}\right\}\right]. \label{GCA-M-def2}
\end{split}
\end{align}
\end{prop}
\begin{proof}
Applying \eqref{cmax} to \eqref{GCA-C-mut}, \eqref{GCA-G-mut}, and \eqref{GCA-M-def}, we have the claim.
\end{proof}

\subsection{Separation formulas for normalized generalized cluster patterns}
In this subsection, we provide the separation formulas for normalized generalized cluster patterns.

\begin{thm}[Separation formulas for normalized generalized cluster patterns] \label{prop:GCA-sep} For a normalized generalized cluster pattern $\bsig=(\Sigma_{t})_{t\in\T_n}$, the following hold.
\begin{align}
y_{i\ell;t}&=\left(\prod_{j=1}^n\prod_{m=1}^{r_j}y_{jm}^{c_{m\ell;t}^{(ji)}}\right)\prod_{j=1}^n \left.F_{j;t}\right|_\PP(\by)^{\ell b_{ji;t}},\label{GCA:y-sep}\\
\hy_{i\ell;t}&=\left(\prod_{j=1}^n\prod_{m=1}^{r_j}\hy_{jm}^{c_{m\ell;t}^{(ji)}}\right)\prod_{j=1}^n F_{j;t}(\hby)^{\ell b_{ji;t}},\label{GCA:hy-sep}\\
x_{i;t}&=\left(\prod_{j=1}^nx_j^{g_{ji;t}}\right)\dfrac{F_{i;t}(\hby)}{F_{i;t}|_\PP(\by)}. \label{GCA:x-sep}
\end{align}
\end{thm}
\begin{proof}
We can prove this by an argument similar to that for ordinary cluster patterns, as follows.
First, we show the formula for $y$-variables. For $t=t_0$, it clearly holds from the initial values of $C$-matrices and $F$-polynomials. Next, we assume that \eqref{GCA:y-sep} holds at $t\in\T_n$ and let $t'$ be $k$-adjacent to $t$. We denote the right-hand side of (\ref{GCA:y-sep}) by $y_{i\ell;t'}$. For $i=k$, we have
\begin{align}
y_{k\ell;t'}&=\left(\prod_{j=1}^n\prod_{m=1}^{r_j}y_{jm}^{c_{m\ell;t'}^{(jk)}}\right)\prod_{j=1}^nF_{j;t'}|_\PP(\by)^{\ell b_{jk;t'}}\notag\\
&=\left(\prod_{j=1}^n\prod_{m=1}^{r_j}y_{jm}^{c_{m,r_k-\ell;t}^{(jk)}-c_{mr_k;t}^{(jk)}}\right)\prod_{j=1}^nF_{j;t}|_\PP(\by)^{-(r_k-(r_k-\ell)) b_{j k;t}}\notag\\
&=\left(\prod_{j=1}^n\prod_{m=1}^{r_j}y_{jm}^{c_{m,r_k-\ell;t}^{(jk)}}\right)\prod_{j=1}^nF_{j;t}|_\PP(\by)^{(r_k-\ell) b_{j k;t}}\times\left(\prod_{j=1}^n\prod_{m=1}^{r_j}y_{jm}^{-c_{mr_k;t}^{(jk)}}\right)\prod_{j=1}^nF_{j;t}|_\PP(\by)^{-r_k b_{jk;t}}\notag\\
&\overset{\eqref{GCA:y-sep}}{=}\dfrac{y_{k,r_k-\ell;t}}{y_{k,r_k;t}},
\end{align}
thus it holds. For $i\neq k$, 
\begin{align}
y_{i\ell;t'}&=\left(\prod_{j=1}^n\prod_{m=1}^{r_j}y_{jm}^{c_{m\ell;t'}^{(ji)}}\right)\prod_{j=1}^nF_{j;t'}|_\PP(\by)^{\ell b_{ji;t'}}\notag\\
&=\left(\prod_{j=1}^n\prod_{m=1}^{r_j}y_{jm}^{c_{m\ell;t}^{(ji)}+c_{mr_k;t}^{(jk)}[\ell b_{ki;t}]_++[-c_{mr_k;t}^{(jk)}]_+(\ell b_{ki;t})}\right)\notag\\
&\quad\times\prod_{\substack{j=1\\j \neq k}}^nF_{j;t}|_\PP(\by)^{\ell( b_{j i;t}+r_k b_{jk;t}[ b_{ki;t}]_++r_k[- b_{j k;t}]_+ b_{ki;t})}\notag\\
&\quad\times\left(F_{k;t}|_\PP(\by)^{-1}\left(\prod_{j=1}^nF_{j;t}|_\PP(\by)^{[-r_k b_{j k;t}]_+}\right)\left(\prod_{j=1}^n\prod_{m=1}^{r_j}y_{jm}^{[-c_{mr_k;t}^{(jk)}]_+}\right)\right.\notag\\
&\quad\times\left.\left[1\op\bigoplus_{w=1}^{r_k}\left\{\left(\prod_{j=1}^n\prod_{m=1}^{r_j}y_{jm}^{c_{mw;t}^{(jk)}}\right)\prod_{j=1}^nF_{j;t}|_\PP(\by)^{w b_{j k;t}}\right\}\right]\right)^{-\ell b_{ki;t}}\notag\\
&\overset{\eqref{GCA:y-sep}}{=}y_{i\ell;t}y_{kr_k;t}^{[\ell b_{ki;t}]_+}(1\op y_{k1;t}\op\cdots\op y_{kr_k;t})^{-\ell b_{ki;t}}.
\end{align}
Therefore, \eqref{GCA:y-sep} also holds for $i\neq k$. The formula for $\hy$-variables holds similarly to that for $y$-variables. Next, we show the formula for $x$-variables. The case of $t=t_0$ is obvious from the initial conditions of $G$-matrices and $F$-polynomials. We assume that \eqref{GCA:x-sep} holds at $t\in\T_n$ and let $t'$ be $k$-adjacent to $t$. When $i\neq k$, it holds since the mutations of $G$-matrices and $F$-polynomials do not change their forms. We denote the right-hand side of (\ref{GCA:x-sep}) by $x_{i;t'}$. For $i=k$, we have
\begin{align}
x_{k;t'}&=\left(\prod_{j=1}^nx_j^{g_{jk;t'}}\right)\dfrac{F_{k;t'}(\hby)}{F_{k;t'}|_\PP(\by)} \notag\\
&=\prod_{j=1}^nx_j^{-g_{jk;t}+\sum_{\ell=1}^ng_{j\ell;t}[-r_kb_{\ell k;t}]_+-\sum_{\ell=1}^n r_\ell b_{j\ell;t_0}\left[-c_{r_\ell r_k;t}^{(\ell k)}\right]_+}\notag\\
&\quad\times\dfrac{F_{k;t}|_\PP(\by)}{F_{k;t}(\hby)}\cdot\dfrac{\prod_{j=1}^nF_{j;t}(\hby)^{[-r_kb_{jk;t}]_+}}{\prod_{j=1}^nF_{j;t}|_\PP(\by)^{[-r_kb_{jk;t}]_+}}\notag\\
&\quad\times\dfrac{\prod_{\ell=1}^n\prod_{m=1}^{r_\ell}\hy_{\ell m}^{[-c_{mr_k;t}^{(\ell k)}]_+}}{\prod_{\ell=1}^n\prod_{m=1}^{r_\ell}y_{\ell m}^{[-c_{mr_k;t}^{(\ell k)}]_+}}\notag\\
&\quad\times\dfrac{1+\sum_{w=1}^{r_k}\left\{\left(\prod_{j=1}^n\prod_{m=1}^{r_j}\hy_{jm}^{c_{mw;t}^{(jk)}}\right)\prod_{j=1}^nF_{j;t}(\hby)^{w b_{jk;t}}\right\}}{1\op\bigoplus_{w=1}^{r_k}\left\{\left(\prod_{j=1}^n\prod_{m=1}^{r_j}y_{jm}^{c_{mw;t}^{(jk)}}\right)\prod_{j=1}^nF_{j;t}|_\PP(\by)^{w b_{jk;t}}\right\}}\notag\\
&\overset{\eqref{GCA:x-sep}}{=}x_{k;t}^{-1}\displaystyle\left(\prod_{j=1}^nx_{j;t}^{[- b_{jk;t}]_+}\right)^{r_k}\dfrac{1+\hy_{k1;t}+\dots+\hy_{kr_k;t}}{1\op y_{k1;t}\op\dots\op y_{kr_k;t}}.
\end{align}
In the last equality, we used \eqref{GCA:y-sep} and \eqref{GCA:hy-sep}. 
\end{proof}

\subsection{Examples}
We calculate the $x$-, $y$-variables, the $C$-, $G$-matrices, and the $F$-polynomials of type $B_2$. Furthermore, we verify the separation formulas and Proposition \ref{last row of C} (2) in these cases.

\begin{ex}[Type $B_2$]\label{B_2 type}  The following is the simplest example in generalized cluster patterns. 
Let $\br=(r_1,r_2)=(2,1), \hat{B}_0=\begin{pmatrix}
0 & -1\\
2 & 0
\end{pmatrix},B_0=\begin{pmatrix}
0 & -1\\
1 & 0
\end{pmatrix}$. Also, let any vertex of $\T_2$ be denoted by $t_i$, and assume $t_{2i}$ and $t_{2i+1}$ in $\T_2$ are $1$-adjacent for any $i\in\Z$. We write $x_{i;j}=x_{i;t_j}$ and $y_{i\ell;j}=y_{i\ell;t_j}$. In this case, expressing the $y$-variables and $x$-variables in terms of the initial $y$-variables $\by=\by_0$, initial $x$-variables $\bx=\bx_0$, and the initial $\hy$-variables $\hby=\hby_0$, we obtain the following:
\begin{align*}
&\begin{cases}
y_{11;1}=y_{11}y_{12}^{-1},\\
y_{12;1}=y_{12}^{-1},\\
y_{21;1}=y_{21}(1\op y_{11}\op y_{12}),
\end{cases}\quad\begin{cases}
y_{11;2}=y_{11}y_{12}^{-1}(1\op y_{21}\op y_{11}y_{21}\op y_{12}y_{21}),\\
y_{12;2}=y_{12}^{-1}(1\op y_{21}\op y_{11}y_{21}\op y_{12}y_{21})^2,\\
y_{21;2}=y_{21}^{-1}(1\op y_{11}\op y_{12})^{-1},
\end{cases}\\
&\begin{cases}
y_{11;3}=y_{11}(1\op y_{21}\op y_{11}y_{21}\op y_{12}y_{21})^{-1},\\
y_{12;3}=y_{12}(1\op y_{21}\op y_{11}y_{21}\op y_{12}y_{21})^{-2},\\
y_{21;3}=y_{12}^{-1}y_{21}^{-1}(1\op 2y_{21}\op y_{11}y_{21} \op y_{21}^2 \op y_{11}y_{21}^2 \op y_{12}y_{21}^2),
\end{cases}\\
&\begin{cases}
y_{11;4}=y_{11}y_{12}^{-1}y_{21}^{-1}(1\op y_{21}),\\
y_{12;4}=y_{12}^{-1}y_{21}^{-2}(1\op y_{21})^2,\\
y_{21;4}=y_{12}y_{21}(1\op 2y_{21}\op y_{11}y_{21} \op y_{21}^2 \op y_{11}y_{21}^2 \op y_{12}y_{21}^2)^{-1},
\end{cases}\\
&\begin{cases}
y_{11;5}=y_{11}y_{21}(1\op y_{21})^{-1},\\
y_{12;5}=y_{12}y_{21}^2(1\op y_{21})^{-2},\\
y_{21;5}=y_{21}^{-1},
\end{cases}\quad
\begin{cases}
y_{11;6}=y_{11},\\
y_{12;6}=y_{12},\\
y_{21;6}=y_{21},
\end{cases}\\
&\begin{cases}
x_{1;1}=x_1^{-1}\dfrac{1+\hy_{11}+\hy_{12}}{1\op y_{11} \op y_{12}},\\
x_{2;1}=x_2,
\end{cases}\quad \begin{cases}
x_{1;2}=x_1^{-1}\dfrac{1+\hy_{11}+\hy_{12}}{1\op y_{11} \op y_{12}},\vspace{0.5em}\\
x_{2;2}=x_2^{-1}\dfrac{1+\hy_{21}+\hy_{11}\hy_{21}+\hy_{12}\hy_{21}}{1\op y_{21}\op y_{11}y_{21}\op y_{12}y_{21}},
\end{cases}\\
&\begin{cases}
x_{1;3}=x_1x_2^{-2}\dfrac{1+2\hy_{21}+\hy_{11}\hy_{21}+\hy_{21}^2+\hy_{11}\hy_{21}^2+\hy_{12}\hy_{21}^2}{1\op 2y_{21}\op y_{11}y_{21}\op y_{21}^2\op y_{11}y_{21}^2\op y_{12}y_{21}^2}, \vspace{0.5em}\\
x_{2;3}=x_2^{-1}\dfrac{1+\hy_{21}+\hy_{11}\hy_{21}+\hy_{12}\hy_{21}}{1\op y_{21}\op y_{11}y_{21}\op y_{12}y_{21}},
\end{cases}\\
&\begin{cases}
x_{1;4}=x_1x_2^{-2}\dfrac{1+2\hy_{21}+\hy_{11}\hy_{21}+\hy_{21}^2+\hy_{11}\hy_{21}^2+\hy_{12}\hy_{21}^2}{1\op 2y_{21}\op y_{11}y_{21}\op y_{21}^2\op y_{11}y_{21}^2\op y_{12}y_{21}^2}, \vspace{0.5em}\\
x_{2;4}=x_1x_2^{-1}\dfrac{1+\hy_{21}}{1\op y_{21}},
\end{cases}\\
&\begin{cases}
x_{1;5}=x_1,\\
x_{2;5}=x_{2;4}=x_1x_2^{-1}\dfrac{1+\hy_{21}}{1\op y_{21}},
\end{cases} \quad \begin{cases}
x_{1;6}=x_1,\\
x_{2;6}=x_2.
\end{cases}
\end{align*}
On the other hand, the $C$-, $G$-matrices, and $F$-polynomials are obtained as follows:
\begin{align*}
&C_1 = \left(\begin{array}{cc|c} 1 & 0 & 0 \\ -1 & -1 & 0 \\ \hline 0 & 0 & 1 \end{array}\right),\quad G_1 = \begin{pmatrix} -1 & 0 \\ 0 & 1 \end{pmatrix},\quad F_{1;1}(\by) = 1+y_{11}+y_{12},\quad F_{2;1}(\by)=1,\\
&C_2 = \left(\begin{array}{cc|c} 1 & 0 & 0 \\ -1 & -1 & 0 \\ \hline 0 & 0 & -1 \end{array}\right),\quad G_2 = \begin{pmatrix} -1 & 0 \\ 0 & -1 \end{pmatrix},\\
&F_{1;2}(\by)=1+y_{11}+y_{12},\quad F_{2;2}(\by)=1+y_{21}+y_{12}y_{21}+y_{11}y_{21},\\
&C_3 = \left(\begin{array}{cc|c} 1 & 0 & 0 \\ 0 & 1 & -1 \\ \hline 0 & 0 & -1 \end{array}\right), \quad G_3 = \begin{pmatrix} 1 & 0 \\ -2 & -1 \end{pmatrix},\\
&F_{1;3}(\by)=1+2y_{21}+y_{21}^2+y_{12}y_{21}^2+y_{11}(y_{21}+y_{21}^2), \quad F_{2;3}(\by) =1+y_{21}+y_{12}y_{21}+y_{11}y_{21},\\
&C_4 = \left(\begin{array}{cc|c} 1 & 0 & 0 \\ -1 & -1 & 1 \\ \hline -1 & -2 & 1 \end{array}\right), \quad G_4 = \begin{pmatrix} 1 & 1 \\ -2 & -1 \end{pmatrix},\\
&F_{1;4}(\by)=1+2y_{21}+y_{21}^2+y_{12}y_{21}^2+y_{11}(y_{21}+y_{21}^2),\quad F_{2;4}(\by) = 1 + y_{21},\\
&C_5 = \left(\begin{array}{cc|c} 1 & 0 & 0 \\ 0 & 1 & 0 \\ \hline 1 & 2 & -1 \end{array}\right), \quad G_5 = \begin{pmatrix} 1 & 1 \\ 0 & -1 \end{pmatrix},\quad F_{1;5}(\by) = 1,\quad F_{2;5}(\by) = 1 + y_{21},\\
&C_6 = \left(\begin{array}{cc|c} 1 & 0 & 0 \\ 0 & 1 & 0 \\ \hline 0 & 0 & 1 \end{array}\right), \quad G_6 = \begin{pmatrix} 1 & 0 \\ 0 & 1 \end{pmatrix},\quad F_{1;6}(\by)=F_{2;6}(\by)=1.
\end{align*}
Comparing with the $x$- and $y$-variables, we confirm that the separation formulas hold.
\end{ex}

\begin{ex} Let $\br=(2,3)$, $\hat{B}_0=\begin{pmatrix}
0 & -3\\
2 & 0
\end{pmatrix}$, $B_0=\begin{pmatrix}
0 & -1\\
1 & 0
\end{pmatrix}$, and the labels of the vertices and edges of $\T_2$ be the same as in Example \ref{B_2 type}. We have $C_k$ for $k=1,\dots,6$ as follows.

\begin{align*}
C_1 &= \left(\begin{array}{cc|ccc}
1 & 0 & 0 & 0 & 0 \\
-1 & -1 & 0 & 0 & 0 \\ \hline
0 & 0 & 1 & 0 & 0 \\
0 & 0 & 0 & 1 & 0 \\
0 & 0 & 0 & 0 & 1
\end{array}\right), \hh C_2 = \left(\begin{array}{cc|ccc}
1 & 0 & 0 & 0 & 0 \\
-1 & -1 & 0 & 0 & 0 \\ \hline
0 & 0 & 0 & 1 & 0 \\
0 & 0 & 1 & 0 & 0 \\
0 & 0 & -1 & -1 & -1
\end{array}\right), \hh C_3 = \left(\begin{array}{cc|ccc}
1 & 0 & 0 & 0 & 0 \\
0 & 1 & -1 & -2 & -3 \\ \hline
0 & 0 & 0 & 1 & 0 \\
0 & 0 & 1 & 0 & 0 \\
0 & 0 & -1 & -1 & -1
\end{array}\right), \\
C_4 &= \left(\begin{array}{cc|ccc}
1 & 0 & 0 & 0 & 0 \\
-3 & -5 & 1 & 2 & 3 \\ \hline
0 & 0 & 1 & 0 & 0 \\
0 & 0 & 0 & 1 & 0 \\
-1 & -2 & 0 & 0 & 1
\end{array}\right),\hh C_5=\left(\begin{array}{cc|ccc}
1 & 0 & 0 & 0 & 0 \\
2 & 5 & -4 & -8 & -12 \\ \hline
0 & 0 & 1 & 0 & 0 \\
0 & 0 & 0 & 1 & 0 \\
1 & 2 & -2 & -4 & -5
\end{array}\right), \hh C_6= \left(\begin{array}{cc|ccc}
1 & 0 & 0 & 0 & 0 \\
-10 & -19 & 4 & 8 & 12 \\ \hline
0 & 0 & 0 & 1 & 0 \\
0 & 0 & 1 & 0 & 0 \\
-4 & -8 & 1 & 3 & 5
\end{array}\right).
\end{align*}
It can be seen that the claim regarding the block structure in Proposition \ref{prop:non-nega} (1) holds. In particular, each time a mutation in direction 2 is applied, the upper-left $2\times2$ submatrix of the $(2,2)$-block alternates between $I_2$ and $\tilde{I}_2$. 
\end{ex}

\subsection{Separation formulas for \texorpdfstring{$\CP$}{P}-normalized generalized cluster patterns}
In the following, we consider non-normalized generalized cluster patterns. First, we define $J$-functions as in the case of ordinary cluster patterns.

\begin{defi} For any $\CP$-normalized generalized cluster  pattern, we define the \emph{generalized $J$-functions} by the following recurrence relations:
\begin{align}
J_{i;t_0}(\by,\CP)&=1,\\
J_{i;t'}(\by,\CP)&=\begin{cases}
\dfrac{N_{k;t}(\by,\CP)}{p_{k;t}J_{k;t}(\by,\CP)} & \text{$(i=k)$},\\
J_{i;t}(\by,\CP) & \text{$(i\neq k)$},
\end{cases}
\end{align}
where
\begin{align*}
\displaystyle N_{k;t}(\by,\CP)&=\left(\prod_{j=1}^nJ_{j;t}(\by,\CP)^{[-r_k b_{jk;t}]_+}\right)\left(\prod_{j=1}^n\prod_{m=1}^{r_j}y_{jm}^{[-c_{mr_k;t}^{(jk)}]_+}\right)\\
&\quad\times\left[1+
\sum_{\ell=1}^{r_k}\left\{\left(\prod_{j=1}^n\prod_{m=1}^{r_j}y_{jm}^{c_{m\ell;t}^{(jk)}}\right)\prod_{j=1}^nJ_{j;t}(\by,\CP)^{\ell b_{jk;t}}\right\}\right].
\end{align*}
\end{defi}\vv

$M_{k;t}(\by)$ is obtained from $N_{k;t}(\by,\CP)$ by the specialization setting all $p_{j;t}$ to $1$: $M_{k;t}(\by)=N_{k;t}(\by,\one)$. Thus, the $F$-polynomials are obtained from the $J$-functions by specializing all $p_{j;t}$ to $1$: $F_{i;t}(\by)=J_{i;t}(\by,\one)$. Unlike the $F$-polynomials, the $J$-functions are rational functions in $\by$ and $\CP$.

Using the $J$-functions, the following separation formulas hold.

\begin{thm} [Separation formulas for $\CP$-normalized generalized cluster patterns] \label{GCA-sep-p} For a $\CP$-normalized generalized cluster pattern, the following separation formulas hold, respectively.
\begin{align}
y_{i\ell;t}&=\left(\prod_{j=1}^n\prod_{m=1}^{r_j}y_{jm}^{c_{m\ell;t}^{(ji)}}\right)\prod_{j=1}^n J_{j;t}|_\PP(\by,\CP)^{\ell b_{ji;t}}\label{GCA-y-sep-p},\\
\hy_{i\ell;t}&=\left(\prod_{j=1}^n\prod_{m=1}^{r_j}\hy_{jm}^{c_{m\ell;t}^{(ji)}}\right)\prod_{j=1}^n F_{j;t}(\hby)^{\ell b_{ji;t}}\label{GCA-hy-sep-p},\\
x_{i;t}&=\left(\prod_{j=1}^nx_j^{g_{ji;t}}\right)\dfrac{F_{i;t}(\hby)}{J_{i;t}|_\PP(\by,\CP)}\label{GCA-x-sep-p}.
\end{align}
\end{thm}
\begin{proof}
This is shown similarly to the case of normalized generalized cluster patterns.
\end{proof}

\section{Further properties}\label{Further}
In this section, we observe some properties regarding normalized generalized cluster patterns. In particular, we look at their correspondence with ordinary cluster patterns and results parallel to the dualities in ordinary cluster patterns.

\subsection{Correspondence between $C,G,F$ for ordinary cluster patterns and $C,G,F$ for generalized cluster patterns}\label{subsec:corres}
In the following, we examine the correspondence between the \texorpdfstring{$C$}{C}-, \texorpdfstring{$G$}{G}-matrices, and \texorpdfstring{$F$}{F}-polynomials in ordinary cluster patterns and those in generalized cluster patterns.\vv

Let $B_t$ be an exchange matrix, and let $\hat{B}_t=B_tR$. Let $C_t,G_t$, and $F_{i;t}(\by)$ be the $C$-, $G$-matrices, and the $F$-polynomials in generalized cluster patterns with the initial exchange matrix $B_{t_0}$ respectively, and let $\hat{C}_t,\hat{G}_t$, and $\hat{F}_{i;t}(\by)$ be the $C$-, $G$-matrices, and $F$-polynomials in ordinary cluster patterns with the same initial vertex and the initial exchange matrix $\hat{B}_{t_0}$ respectively. Furthermore, let $C^*_t$ be the $n \times n$ matrix whose $(i,j)$-entry is the $(r_i,r_j)$-entry of the $(i,j)$-block of $C_t$. We also define $F^*_{i;t}(\by)$ to be the expression obtained from $F_{i;t}(\by)$ by setting all initial $y$-variables except $y_{ir_i}$ to $0$, and replacing each $y_{ir_i}$ with $y_i$. This is possible because $\by$ is a formal variable in $\PP$. Then, the following hold.

\begin{thm}\label{thm:coincide} \textup{(1)} For any $t \in \T_n$, we have $C^*_t = \hat{C}_t$ and $G_t = \hat{G}_t$.\\
\textup{(2)} For any $i\in\{1,\dots,n\}$ and $t\in\T_n$, $F^*_{i;t}(\by) = \hat{F}_{i;t}(\by)$ holds.
\end{thm}
\begin{proof}
We have the claim (1) by Proposition \ref{prop:new-mut}. Below, we show (2). From Proposition \ref{Block-str}, $M_{k;t}(\by)$ in \eqref{GCA-M-def2} can be written as follows:
\begin{align}
\begin{split}
\displaystyle M_{k;t}(\by)&=\left(\prod_{j=1}^nF_{j;t}(\by)^{[-\hat{b}_{jk;t}]_+}\right)\left(\prod_{j=1}^n y_{jr_j}^{[-c_{r_jr_k;t}^{(jk)}]_+}\right)\\
&\quad \times \left[1+
\sum_{\ell=1}^{r_k-1}\left\{y_{k\ell}^{c_{\ell\ell;t}^{(kk)}}y_{k,r_k-\ell}^{c_{r_k-\ell,\ell;t}^{(kk)}}\left(\prod_{j=1}^n y_{jr_j}^{c_{r_j\ell;t}^{(jk)}}\right)\prod_{j=1}^nF_{j;t}(\by)^{\ell b_{jk;t}}\right\}+\left(\prod_{j=1}^n y_{jr_j}^{c_{r_jr_k;t}^{(jk)}}\right)\prod_{j=1}^nF_{j;t}(\by)^{\hat{b}_{jk;t}}\right].\label{GCA:M-3}
\end{split}
\end{align}
In view of Proposition \ref{Block-str}, $y_{k\ell}^{c_{\ell\ell;t}^{(kk)}}y_{k,r_k-\ell}^{c_{r_k-\ell,\ell;t}^{(kk)}}$ is either $y_{k\ell}$ or $y_{k,r_k-\ell}$. Therefore, it vanishes under the specialization of $\by$. Then, using the claim (1), we observe that the mutation formulas for $F_{i;t}^\ast(\by)$ coincide with that for the ordinary cluster patterns.
\end{proof}

\begin{cor}\label{const=1} The constant term of each generalized $F$-polynomial $F_{i;t}(\by)$ is $1$. 
\end{cor}
\begin{proof}
By Theorem \ref{thm:coincide}, the constant term of a generalized $F$-polynomial $F_{i;t}(\by)$ coincides with the constant term of the corresponding $F$-polynomial $\hat{F}_{i;t}(\by)$ in an ordinary cluster pattern, and the constant term of $\hat{F}_{i;t}(\by)$ is $1$ by \cite{GHKK18}.
\end{proof}

Next, we consider the determinants of generalized $C$-matrices. Let $\tilde{\Lambda}_n=\{i\in\{1,\dots,n\}\mid r_i-1\equiv 2,3 \pmod 4\}$ and $\sigma_t=\prod_{i\in \tilde{\Lambda}_n}\sigma_{i;t}$. Then, the following holds.

\begin{prop}[Unimodularity] We have
\begin{align}
\det(C_t)=\sigma_t \det(\hat{C}_t).
\end{align}
In particular, $C$-matrices are unimodular. 
\end{prop}

\begin{proof}
First, we note that $\tilde{\Lambda}_n=\{i\in\{1,\dots,n\}\mid \det(\tilde{I}_{r_i-1})=-1\}$. By Proposition~\ref{Block-str}, $C_t$ is expressed in the following form:
\begin{align*}
C_t=
\left(\begin{array}{cc|c|c|c}
I_{r_1-1} \text{ or } \tilde{I}_{r_1-1} & \bze & O & \cdots & O\\ 
\ast & c_{r_1r_1;t}^{(11)} & \ast & & \ast\\ \hline
\multicolumn{2}{c|}{C_{21;t}} & C_{22;t} & \cdots & C_{2n;t}\\ \hline
\multicolumn{2}{c|}{\vdots} & \vdots & \ddots & \vdots\\ \hline
\multicolumn{2}{c|}{C_{n1;t}} & C_{n2;t} & \cdots & C_{nn;t}
\end{array}\right).
\end{align*} 
By cofactor expansion of $\det(C_t)$ along the first row, we have
\begin{align*}
\det(C_t)=(c_{11;t}^{(11)}+(-1)^{r_1}c_{1,r_1-1;t}^{(11)})
\det \left(\begin{array}{cc|c|c|c}
I_{r_1-2} \text{ or } \tilde{I}_{r_1-2} & \bze & O & \cdots & O\\ 
\ast & c_{r_1r_1;t}^{(11)} & \ast & & \ast\\ \hline
\multicolumn{2}{c|}{C_{21;t}'} & C_{22;t} & \cdots & C_{2n;t}\\ \hline
\multicolumn{2}{c|}{\vdots} & \vdots & \ddots & \vdots\\ \hline
\multicolumn{2}{c|}{C_{n1;t}'} & C_{n2;t} & \cdots & C_{nn;t}
\end{array}\right),
\end{align*} 
where $C_{i1;t}'$ is the $r_i \times (r_1-1)$ matrix obtained from $C_{i1;t}$ by deleting either the first column or the $(r_1-1)$-th column.
Repeating this procedure, we obtain
\begin{align*}
\det(C_t)&=(c_{11;t}^{(11)}+(-1)^{r_1}c_{1,r_1-1;t}^{(11)})
\det \left(\begin{array}{cc|c|c|c}
I_{r_1-2} \text{ or } \tilde{I}_{r_1-2} & \bze & O & \cdots & O\\ 
\ast & c_{r_1r_1;t}^{(11)} & \ast & & \ast\\ \hline
\multicolumn{2}{c|}{C_{21;t}'} & C_{22;t} & \cdots & C_{2n;t}\\ \hline
\multicolumn{2}{c|}{\vdots} & \vdots & \ddots & \vdots\\ \hline
\multicolumn{2}{c|}{C_{n1;t}'} & C_{n2;t} & \cdots & C_{nn;t}
\end{array}\right)\\
&\;\:\vdots\\
&=\det (I_{r_1-1} \text{ or } \tilde{I}_{r_1-1})
\det \left(\begin{array}{c|c|c|c}
c_{r_1r_1;t}^{(11)} & \ast & & \ast\\ \hline
\bc_{r_1;t}^{(21)} & C_{22;t} & \cdots & C_{2n;t}\\ \hline
\vdots & \vdots & \ddots & \vdots\\ \hline
\bc_{r_1;t}^{(n1)} & C_{n2;t} & \cdots & C_{nn;t}
\end{array}\right).
\end{align*} 
Applying this procedure to the remaining $i$-th row blocks for $i=2,\dots,n$, we obtain $\det(C_t)=\sigma_t \det(\hat{C}_t)$. Furthermore, since $\det(\hat{C}_t)=\pm1$ by \cite{NZ12}, it follows that $\det(C_t)=\pm1$.
\end{proof}

\begin{ex}  The $C$-, $G$-matrices and $F$-polynomials in an ordinary cluster pattern with the initial exchange matrix $\hat{B}_0$ from Example \ref{B_2 type} $\hat{C}_t, \hat{G}_t, \hat{F}_{i;t}(\by)$ are as follows:
\begin{align*}
&\hat{C}_1 =\begin{pmatrix}-1 & 0 \\ 0 & 1 \end{pmatrix},\quad \hat{G}_1 = \begin{pmatrix} -1 & 0 \\ 0 & 1 \end{pmatrix},\quad \hat{F}_{1;1}(\by) = 1+y_{1},\quad \hat{F}_{2;1}(\by)=1,\\
&\hat{C}_2 =\begin{pmatrix}-1 & 0 \\ 0 & -1 \end{pmatrix},\quad \hat{G}_2 = \begin{pmatrix} -1 & 0 \\ 0 & -1 \end{pmatrix},\quad \hat{F}_{1;2}(\by)=1+y_{1},\quad \hat{F}_{2;2}(\by)=1+y_{2}+y_{1}y_{2},\\
&\hat{C}_3 =\begin{pmatrix}1 & -1 \\ 0 & -1 \end{pmatrix}, \quad \hat{G}_3 = \begin{pmatrix} 1 & 0 \\ -2 & -1 \end{pmatrix},\\
&\hat{F}_{1;3}(\by)=1+2y_{2}+y_{2}^2+y_{1}y_{2}^2, \quad \hat{F}_{2;3}(\by) = 1+y_{2}+y_{1}y_{2},\\
&\hat{C}_4 =\begin{pmatrix}-1 & 1 \\ -2 & 1 \end{pmatrix}, \quad \hat{G}_4 = \begin{pmatrix} 1 & 1 \\ -2 & -1 \end{pmatrix},\\
&\hat{F}_{1;4}(\by)=1+2y_{2}+y_{2}^2+y_{1}y_{2}^2,\quad \hat{F}_{2;4}(\by) = 1 + y_{2},\\
&\hat{C}_5 =\begin{pmatrix}1 & 0 \\ 2 & -1 \end{pmatrix}, \quad \hat{G}_5 = \begin{pmatrix} 1 & 1 \\ 0 & -1 \end{pmatrix},\quad \hat{F}_{1;5}(\by) = 1,\quad \hat{F}_{2;5}(\by) = 1 + y_{2},\\
&\hat{C}_6 =\begin{pmatrix}1 & 0 \\ 0 & 1 \end{pmatrix}, \quad \hat{G}_6 = \begin{pmatrix} 1 & 0 \\ 0 & 1 \end{pmatrix},\quad \hat{F}_{1;6}(\by)=\hat{F}_{2;6}(\by)=1.
\end{align*}
By comparing them with $C_t, G_t, F_{i;t}(\by)$ in Example \ref{B_2 type}, we observe that Theorem \ref{thm:coincide} holds.
\end{ex}

\subsection{Dualities of \texorpdfstring{$C$}{C}- and \texorpdfstring{$G$}{G}-matrices}\label{subsec:dual}

We examine the dualities between generalized $C$- and $G$-matrices as in the ordinary cases. To describe the formula in a parallel way to the ordinary case, we introduce the following matrices.
\begin{defi}\label{def:Gamma-mat}
We define the \emph{$\gamma$-vector} $\boldsymbol{\gamma}_{j;t}$ as follows:
\begin{align}
\boldsymbol{\gamma}_{j;t}=r_j^{-1}\begin{pmatrix}
r_1c_{r_1r_j;t}^{(1j)}\\
\vdots\\
r_nc_{r_nr_j;t}^{(nj)}
\end{pmatrix}. \label{eq:gamma-vector}
\end{align}
Each entries of $\boldsymbol{\gamma}_{j;t}$ is integer by \eqref{eq:last row of C} for $\ell=1$. The $n\times n$ matrix $\Gamma_t=(\gamma_{ij;t})$ whose $j$-th column vector is the $\gamma$-vector $\boldsymbol{\gamma}_{j;t}$ is called the $\Gamma$-matrix.
\end{defi}

By the definition of $\Gamma$-matrices, we have the following relation to $C$-matrices $C_t^\ast$ in Theorem \ref{thm:coincide}:
\begin{align}
\Gamma_t=R C_t^\ast R^{-1},\quad \gamma_{ij;t}=r_ir_j^{-1}c_{r_ir_j;t}^{(ij)}. \label{C and gamma}
\end{align}

\begin{prop}[First duality] \label{prop:first-dual} The following equality holds.
\begin{align}
B_{t_0}RC_t^\ast=G_tB_tR.\label{first}
\end{align}
\end{prop}
\begin{proof}
By \cite[(6.14)]{CA4}, we have $\hat{B}_{t_0}\hat{C}_t=G_t\hat{B}_t$. Then, the claim holds by Theorem \ref{thm:coincide} (1).
\end{proof}

By \eqref{C and gamma}, equation \eqref{first} is expressed in a parallel form to \cite[(6.14)]{CA4} as follows:
\begin{align}
B_{t_0}\Gamma_t=G_tB_t. \label{form-first}
\end{align}

\begin{ex} Let $\br=(2,3)$,\quad$\hat{B}_0=\begin{pmatrix}
0 & -3\\
2 & 0
\end{pmatrix}$,\quad$B_0=\begin{pmatrix}
0 & -1\\
1 & 0
\end{pmatrix}$. Then, we have
\begin{align*}
C_8 = \left(\begin{array}{cc|ccc}
1 & 0 & 0 & 0 & 0 \\
-36 & -71 & 15 & 30 & 45 \\ \hline
0 & 0 & 1 & 0 & 0 \\
0 & 0 & 0 & 1 & 0 \\
-15 & -30 & 6 & 12 & 19
\end{array}\right),\quad 
\Gamma_8=\begin{pmatrix}
-71 & 30\\
-45 & 19
\end{pmatrix},\quad
G_8 = \begin{pmatrix}
19 & 45 \\
-30 & -71
\end{pmatrix}.
\end{align*}
Thus, we obtain $B_{0}\Gamma_8=G_8B_8=\begin{pmatrix}
45 & -19 \\
-71 & 30
\end{pmatrix}$. Also, let $\hat{C}_n$ and $\hat{G}_n$ be the $C$- and $G$-matrix in ordinary cluster patterns with the initial exchange matrix $\hat{B}_0$, respectively. Then, we have $\hat{C}_8=\begin{pmatrix}
-71 & 45\\
-30 & 19
\end{pmatrix}, \hat{G}_8=\begin{pmatrix}
19 & 45 \\
-30 & -71
\end{pmatrix}$, and we confirm that the $G$-matrices coincide, and the $(r_i,r_j)$-entry of $C_8$ is the $(i,j)$-entry of $\hat{C}_8$. \end{ex}

\begin{prop} \label{prop:gamma-mut} Let $t'$ be $k$-adjacent to $t$. Then, the matrix $\Gamma_t=(\gamma_{ij;t})$ mutates in direction $k$ as follows:
\begin{align}
\gamma_{ij;t'}=\begin{cases}
-\gamma_{ik;t} & \text{$(j=k)$},\\
\gamma_{ij;t}+r_k\gamma_{ik;t}[ b_{kj;t}]_++r_k[-\gamma_{ik;t}]_+ b_{kj;t} & \text{$(j\neq k)$}.
\end{cases}\label{gamma-mut}
\end{align}
\end{prop}
\begin{proof}
By \eqref{C and gamma}, we have $\gamma_{ij;t}=r_j^{-1}r_ic_{r_ir_j;t}^{(ij)}$. Then, for $j\neq k$, 
\begin{align*}
\gamma_{ij;t'}&=r_j^{-1}r_i(c_{r_ir_j;t}^{(ij)}+c_{r_ir_k;t}^{(ik)}[r_j b_{kj;t}]_++[-c_{r_ir_k;t}^{(ik)}]_+(r_j b_{kj;t}))\\
&=\gamma_{ij;t}+r_k\cdot r_k^{-1}r_ic_{r_ir_k;t}^{(ik)}[ b_{kj;t}]_++r_k\cdot r_k^{-1}r_i[-c_{r_ir_k;t}^{(ik)}]_+ b_{kj;t}\\
&=\gamma_{ij;t}+r_k[ b_{kj;t}]_+\gamma_{ik;t}+r_k[-\gamma_{ik;t}]_+ b_{kj;t}.
\end{align*}
Also, since $c_{r_ir_k;t'}^{(ik)}=-c_{r_ir_k;t}^{(ik)}$, $\gamma_{ik;t'}=-\gamma_{ik;t}$ holds.
\end{proof}

The following lemma gives $\vep$-expressions for $C$-, $\Gamma$-, and $G$-matrices in generalized cluster patterns in parallel with those for ordinary cluster patterns \cite[(2.4) and (2.5)]{NZ12}.

\begin{lem}\label{GCACGep} The following three expressions do not depend on the choice of $\vep\in\{\pm1\}$.
\begin{align}
c_{m\ell;t}^{(ij)}+c_{mr_k;t}^{(ik)}[\vep\ell b_{kj;t}]_++[-\vep c_{mr_k;t}^{(ik)}]_+(\ell b_{kj;t}),\label{GCA-C-mut-ep}\\
\gamma_{ij;t}+r_k[\vep  b_{kj;t}]_+\gamma_{ik;t}+r_k[-\vep \gamma_{ik;t}]_+ b_{kj;t},\label{gamma-mut-ep}\\
-g_{ik;t}+\displaystyle\sum_{\ell=1}^ng_{i\ell;t}[-\vep \hat{b}_{\ell k;t}]_+-\displaystyle\sum_{\ell=1}^n \hat{b}_{i\ell;t_0}\left[-\vep c_{r_\ell r_k;t}^{(\ell k)}\right]_+. \label{GCA-G-mut-ep}
\end{align}
\end{lem}

\begin{proof}
This is obvious for \eqref{GCA-C-mut-ep} and \eqref{gamma-mut-ep}. The expression \eqref{GCA-G-mut-ep} does not depend on $\vep\in\{\pm1\}$ from Theorem \ref{thm:coincide} and \cite[(2.5)]{NZ12}.
\end{proof}

\begin{defi} [Tropical signs]
By \cite{GHKK18}, the $C$-matrices $\hat{C}_t$ in ordinary cluster patterns have the column sign-coherence. Then, we define the \emph{tropical sign} $\vep_{k;t}\in\{\pm1\}$ as the sign of the $c$-vector $\hat{\bc}_{i;t}$.
\end{defi}

For a matrix $A$, we denote by $A^{k\bullet}$ (resp. $A^{\bullet k}$) the matrix obtained by replacing all entries of $A$ except for those in the $k$-th row (resp. $k$-th column) with $0$. In addition, we define the $n\times n$ matrix $J_k$ obtained by replacing the $(k,k)$-entry of $n\times n$ identity matrix by $-1$. 

\begin{prop}\label{mut_Mat} The following holds.
\begin{align}
C_{t'}^\ast&=C_t^\ast J_k+C_t^\ast[\vep_{k;t} B_t^{k\bullet}R]_+,\\
\Gamma_{t'}&=\Gamma_tJ_k+\Gamma_t[\vep_{k;t}R B_t^{k\bullet}]_+,\\
G_{t'}&=G_tJ_k+G_t[-\vep_{k;t}B_t^{\bullet k}R]_+.
\end{align}
\end{prop}
\begin{proof}
By Lemma \ref{GCACGep}, we have the following:
\begin{align}
c_{r_ir_j;t'}^{(ij)}&=c_{r_ir_j;t}^{(ij)}+c_{r_ir_k;t}^{(ik)}[\vep_{k;t}r_j b_{kj;t}]_+,\label{GCA-C-mut-tropsign}\\
\gamma_{ij;t'}&=\gamma_{ij;t}+r_k\gamma_{ik;t}[\vep_{k;t} b_{kj;t}]_+ \label{gamma-mut-tropsign},\\
g_{ik;t'}&=-g_{ik;t}+\displaystyle\sum_{\ell=1}^ng_{i\ell;t}[-\vep_{k;t} \hat{b}_{\ell k;t}]_+, \label{GCA-G-mut-tropsign}
\end{align}
for $j\neq k$. Therefore, we obtain the claim.
\end{proof}

\begin{prop} [Second duality] \label{Second duality}
The following holds.
\begin{align}
(R D)^{-1}G_t^T(R D)C_t^\ast=I_n.\label{second}
\end{align}
\end{prop}
\begin{proof}
The claim holds from \cite[(3.11)]{NZ12} and Theorem \ref{thm:coincide} (1).

Alternately, we prove it by the induction on $t\in\T_n$ using Proposition \ref{mut_Mat}.
\end{proof}

By using $\Gamma$-matrix, equation (\ref{second}) can be expressed in a form parallel to \cite[(3.11)]{NZ12} as follows:
\begin{align}
D^{-1}G_t^TD\Gamma_t=I_n. \label{form-second}
\end{align}

\subsection{Correspondence with the result of \cite{Nak15,NR16}}\label{corres}
In the following, we provide the correspondence between our $y$-variables, $C$-, $G$-, and $F$-patterns and those in \cite{Nak15,NR16}. \vv

Table \ref{corres-table} shows the correspondence with the notions in \cite{Nak15,NR16}, where ${}^RC_t$ (resp. ${}^RG_t$) is the ordinary $C$-matrix (resp. $G$-matrix) for the $B$-pattern $\{B_tR\}$, and ${}^LC_t$ (resp. ${}^LG_t$) is the ordinary $C$-matrix (resp. $G$-matrix) for the $B$-pattern $\{RB_t\}$. The result in \cite{Nak15,NR16} is for normalized generalized cluster patterns. Therefore, we assume that generalized cluster patterns are normalized below. The formulation of \cite{Nak15,NR16} needs the assumption $p_{ir_i;t}/p_{i0}=\tilde{y}_{i;t}^{r_i}$, that is, the $r_i$-th root of $p_{ir_i;t}/p_{i0}$ in $\PP$.

\begin{table}[h]
    \centering
    \begin{tabular}{c||ccc}
         & This paper & & \cite{Nak15,NR16}\\ \hline
        Assumption  & None & & $\tilde{y}_{i;t}^{r_i}=p_{ir_i;t}/p_{i0;t}$\\
      $y$-variables & $y_{i\ell;t}$ & &  $z_{i,\ell;t}\tilde{y}_{i;t}^\ell$ \\
     $C$-matrices  & $C_t^\ast=\hat{C}_t$& & $R^{-1}C_tR={}^RC_t$\\
      $\Gamma$- ($C$-) matrices  & $\Gamma_t=RC_t^\ast R^{-1}$ & & $C_t={}^LC_t=R({}^RC_t)R^{-1}$\\
      $G$-matrices & $G_t=\hat{G}_t$ & & $G_t={}^RG_t$\\
      $F$-polynomials & $F_{i;t}(\by)$,\quad $\by=(y_{i\ell})$ & & $F_{i;t}(\tilde{\by},\bz)$,\quad $\tilde{\by}=(\tilde{y}_i)$
    \end{tabular}
    \caption{Correspondence with the notions in \cite{Nak15,NR16}}\label{corres-table}
\end{table}

The coefficient semifield $\Tp(\by,\bz)$ of the generalized cluster pattern with principal coefficients in \cite{Nak15,NR16} is generated by $\tilde{\by}=(\tilde{y}_1,\dots,\tilde{y}_n)$ and the data $\bz=(\bz_1,\dots,\bz_n)$, where $\bz_i=(z_{i,0},\dots,z_{i,r_i})$ satisfies $z_{i,0}=z_{i,r_i}=1$. We compare the two formulations of generalized cluster patterns with principal coefficients, whose coefficient semifields are $\Tp(\by)$ and $\Tp(\tilde{\by},\bz)$, respectively. Below, we show that Proposition \ref{Block-str} and \ref{last row of C} are recovered from the exponents of the $y$-variables in \cite{Nak15,NR16}. For $\ell\neq r_i$, by substituting $z_{i,\ell}\tilde{y}_{i}^\ell$ for the initial variables $y_{i\ell}$ in our formulation, we obtain:
\begin{align}
y_{i\ell;t}&=\prod_{j=1}^n\prod_{m=1}^{r_j}y_{jm}^{c_{m\ell;t}^{(ji)}}=\prod_{j=1}^n\prod_{m=1}^{r_j}(z_{j,m}\tilde{y}_j^m)^{c_{m\ell;t}^{(ji)}}.\label{sub_y}
\end{align}

On the other hand, let ${}^RC_t=(c_{ij;t})$ denote the $C$-matrix (in the right-companion cluster algebra) defined in \cite{Nak15,NR16}. Note that ${}^RC_t=\hat{C}_t=C_t^\ast$. Then, $\tilde{y}_{i;t}^{\ell}$ is given by:
\begin{align}
z_{i,\ell;t}\tilde{y}_{i;t}^\ell&=z_{i,\ell;t}\prod_{j=1}^n(\tilde{y}_j^{c_{ji;t}})^\ell=z_{i,\ell;t}\prod_{j=1}^n\tilde{y}_j^{\ell c_{r_jr_i;t}^{(ji)}}.\label{NR_y_power}
\end{align}
Here, by \cite{NR16}, the data $z_{i,\ell;t}$ mutates in direction $k$ as follows:
\begin{align}
z_{i,\ell;t'}=\begin{cases}
z_{k,r_k-\ell;t} & \text{($i=k$)},\\
z_{i,\ell;t} &\text{($i\neq k$)},
\end{cases}\label{z-mut}
\end{align}
where $t'$ is adjacent to $t$. By \eqref{sub_y} and \eqref{NR_y_power}, we consider the following equation.
\begin{align}
\prod_{j=1}^n\prod_{m=1}^{r_j}(z_{j,m}\tilde{y}_j^m)^{c_{m\ell;t}^{(ji)}}=z_{i,\ell;t}\prod_{j=1}^n\tilde{y}_j^{\ell c_{r_jr_i;t}^{(ji)}}.\label{compare}
\end{align}
In the following, we compare the exponents of $y_j$ in them. First, for $j\neq i$, $z_{j,m}^{c_{m\ell;t}^{(ji)}}$ in the left-hand side of \eqref{compare} must be 1. Thus, we obtain the following:
\begin{align}
mc_{m\ell;t}^{(ji)}=\begin{cases}
\ell c_{r_jr_i;t}^{(ji)} &\text{($m=r_j$)},\\
0 & \text{($m\neq r_j$)}.
\end{cases}\label{compare_1}
\end{align}
Next, we suppose $j=i$ and $\sigma_{i;t}=1$. Since $\sigma_{i;t}=1$, we have $z_{i,\ell;t}=z_{i,\ell}$; thus we have $c_{m\ell;t}^{(ii)}=0$ for $m\neq \ell,r_i$ and $c_{\ell\ell;t}^{(ii)}=1$. Therefore, we obtain
\begin{align}
c_{m\ell;t}^{(ii)}=0\quad (m\neq \ell,r_i), \quad c_{\ell\ell;t}^{(ii)}=1, \quad \ell+r_ic_{r_i\ell;t}^{(ii)}=\ell c_{r_ir_i;t}^{(ii)}.\label{compare_2}
\end{align}
Similarly, for $j=i$ and $\sigma_{i;t}=-1$, we have $z_{i,\ell;t}=z_{i,r_i-\ell}$, $c_{m\ell;t}^{(ii)}=0$ for $m\neq r_i-\ell,r_i$, and $c_{r_i-\ell,\ell;t}^{(ii)}=1$. Thus, we obtain
\begin{align}
c_{m\ell;t}^{(ii)}=0\quad (m\neq r_i-\ell,r_i), \quad c_{r_i-\ell,\ell;t}^{(ii)}=1, \quad r_i-\ell+r_ic_{r_i\ell;t}^{(ii)}=\ell c_{r_ir_i;t}^{(ii)}.\label{compare_3}
\end{align}
Therefore, we recover Proposition \ref{Block-str} and \ref{last row of C} from \eqref{compare_1}, \eqref{compare_2}, \eqref{compare_3}.

\end{document}